\documentclass[a4,12pt]{article} 
\usepackage{graphicx}
\RequirePackage[OT1]{fontenc}
\RequirePackage{amsthm,amsmath,natbib}
\RequirePackage[colorlinks,citecolor=blue,urlcolor=blue]{hyperref}

\numberwithin{equation}{section}
\theoremstyle{plain}
\newtheorem{thm}{Theorem}[section]
\newtheorem{pro}{Proposition}[section]

\newtheorem{lem}{Lemma}[section]
\theoremstyle{remark}
\newtheorem{rem}{Remark}[section]
\newcommand\argmin{\mathrm {argmin}}
\newcommand\Var{\mathrm {Var}}
\newcommand\Cov{\mathrm {Cov}}
\newdimen\AAdi%
\newbox\AAbo%
\def\AArm{\fam0 }%\tenrm}%
\def\AAk#1#2{\setbox\AAbo=\hbox{#2}\AAdi=\wd\AAbo\kern#1\AAdi{}}%
\def\AAr#1#2#3{\setbox\AAbo=\hbox{#2}\AAdi=\ht\AAbo\raise#1\AAdi\hbox{#3}}%
\def\BBn{{\AArm I\!N}}%
\def\BBe{{\AArm I\!E}}
\def\BBr{{\AArm I\!R}}
\def\BBp{{\AArm I\!P}}
\def\BBone{{\AArm 1\AAk{-.8}{I}I}}%
\def\genhGCV{{\hat h_{G}}}
\begin{document}
%\begin{frontmatter}    !!!!!!
\title{On bandwidth selection problems in nonparametric  trend estimation under martingale difference errors}
%\runtitle{Bandwidth selection under stationary MDS errors}     !!!!!!

%\begin{aug}
%\author{\fnms{Karim} \snm{ Benhenni}\thanksref{a,e1} \ead[label=e1,mark]{karim.benhenni@univ-grenoble-alpes.fr}}
%\author{\fnms{Didier A.} \snm{Girard}\thanksref{b,e2} \ead[label=e2,mark]{didier.girard@univ-grenoble-alpes.fr}}
%\and
%\author{\fnms{Sana} \snm{Louhichi}\thanksref{a,e3}%
%\ead[label=e3,mark]{sana.louhichi@univ-grenoble-alpes.fr}%
%%\ead[label=u1,url]{https://ljk.imag.fr/}
%}
%
%\address[a]{
%Laboratoire Jean Kuntzmann.
%Universit\'e Grenoble Alpes, 700 Avenue Centrale,
%38401 Saint-Martin-d'H\`eres, France.
%\printead{e1}
%\printead{e3}
%}
%
%\address[b]{
%CNRS, Laboratoire Jean Kuntzmann, 700 Avenue Centrale,
%38401 Saint-Martin-d'H\`eres, France.
%\printead{e2}
%%, \printead{u1}
%}
%
%\runauthor{K. Benhenni et al.}
%
%\affiliation{Laboratoire Jean Kuntzmann (CNRS 5224)}
%
%\end{aug}

%\bigskip
%{\large{Karim Benhenni\footnote{Laboratoire Jean Kuntzmann,
%Universit\'e Grenoble Alpes, 700 Avenue Centrale,
%38401 Saint-Martin-d'H\`eres, France. E.mail: karim.benhenni@univ-grenoble-alpes.fr}, Didier A. Girard\footnote{CNRS, Laboratoire Jean Kuntzmann.
%E.mail: didier.girard@univ-grenoble-alpes.fr},
%Sana Louhichi\footnote{Corresponding author. Laboratoire Jean Kuntzmann, Universit\'e Grenoble Alpes, 700 Avenue Centrale, 38401 Saint-Martin-d'H\`eres, France. E.mail: sana.louhichi@univ-grenoble-alpes.fr.}}
%}
\author{Karim Benhenni\footnote{Laboratoire Jean Kuntzmann,
Universit\'e Grenoble Alpes, 700 Avenue Centrale,
38401 Saint-Martin-d'H\`eres, France. E.mail: karim.benhenni@univ-grenoble-alpes.fr}, Didier A. Girard\footnote{CNRS, Laboratoire Jean Kuntzmann.
E.mail: didier.girard@univ-grenoble-alpes.fr},
Sana Louhichi\footnote{Laboratoire Jean Kuntzmann, Universit\'e Grenoble Alpes, 700 Avenue Centrale, 38401 Saint-Martin-d'H\`eres, France. E.mail: sana.louhichi@univ-grenoble-alpes.fr.}}

%\title{On smoothing parameters selection problems in nonparametric regression models with martingale difference errors}

%\author{Karim Benhenni\footnote{karim.benhenni@univ-grenoble-alpes.fr}, Didier Girard\footnote{didier.girard@univ-grenoble-alpes.fr},
%Sana Louhichi\footnote{sana.louhichi@univ-grenoble-alpes.fr.\,\,\, Laboratoire Jean Kuntzmann (CNRS 5224), Universit\'e Grenoble Alpes, 700 Avenue Centrale,
%38401 Saint-Martin-d'H\`eres, France.}}

%\maketitle  
\maketitle

\begin{abstract}
In this paper, we are interested in the problem of  smoothing parameter selection in
nonparametric curve estimation under dependent errors. We focus on  kernel estimation and the case when the errors
form a general stationary sequence of martingale difference random variables where neither linearity assumption nor``all moments are finite" are required.
We compare the behaviors of the smoothing bandwidths obtained by minimizing
either the unknown average squared error, the theoretical mean average squared error,
a Mallows-type criterion adapted to the dependent case and the family of criteria known as generalized cross validation (GCV) extensions of the Mallows' criterion. We prove that these three minimizers and those based on the  GCV family are first-order equivalent in probability. We give also a normal asymptotic behavior of the
gap between the minimizer of the average square error and that of the  Mallows-type criterion. This is extended to the GCV family.
Finally, we apply our theoretical results to a specific
case of martingale difference
sequence, namely the Auto-Regressive Conditional Heteroscedastic (ARCH(1)) process.
A Monte-carlo simulation study, for this regression model with ARCH(1) process, is conducted.
\end{abstract}

%\begin{keyword}    !!!!!!
%\kwd{\LaTeXe}

%\end{keyword}    !!!!!!

%\end{frontmatter}      !!!!!!
  \tableofcontents

\section{Introduction}

This paper is about  nonparametric regression model (known also as a machine learning function) which is used as a tool to
  describe and to analyse  the trend between
a response variable and one or more explanatory random variables. This subject was studied by several authors since 1964 (\cite{N},\,\cite{W})
and is still relevant, due to the fact that nonparametric regression has a lot
of  applications in  different fields, such as  economics, medicine, biology, physics, environment,
social sciences, $\cdots$, see for instance \cite{Hastie}.

Several estimate of the nonparametric regression function are proposed in the literature such as
kernel smoothing, local polynomial regression, spline-based regression models, and regression trees (see for instance \cite{Hastie}).
In this paper, we are interested in kernel nonparametric  estimations.
These estimate depend on some smoothing parameter $h$ which has to be chosen according to some criteria.
For independent observations, {\textcolor{black}{two popular criteria}}, to select $h$, are known as  the Cross Validation (CV) criterion and its rotation-invariant
version called
 Generalized Cross-Validation (GCV) criterion. The GCV criterion has different variants, see for instance \cite{A},
 \cite{CW}, \cite{SH}, \cite{roce}, \cite{MA}. We refer the reader to  \cite{Hardele} who studied this problem in the case of independent, equally spaced, observations.
They gave, in particular, the behaviors of the minimizers over $h$ of the average squared errors, the mean average squared errors,
the cross-validation score CV or the generalized cross-validation GCV. They also studied the deviation between these selected smoothing parameters.

{\textcolor{black}{In many cases}}, independence of the observations is, however,  not a realistic modeling of observed data.
% since, in practice, {\textcolor{black}{they are}} often correlated.
Autoregressive models,
autoregressive conditional heteroscedasticity models, Markov chains are examples of dependent models (see for instance \cite{DLLY}).
We focus, in this paper,  on the case of  kernel nonparametric models with particular dependent errors, more precisely, the case when
the errors form a stationary {\textcolor{black}{martingale difference sequence (MDS, in short)}}.
They  are, essentially,  two reasons that motivated us to restrict our  study of dependence to  the case of stationary MDS.
\begin{itemize}
\item[$\bullet$] The first reason is that,  studying MDS  is a promising step for studying the  general case of stationary dependent errors. In fact,
MDS plays an important role in establishing
the results for arbitrary stationary sequences, see for instance  \cite{PUW} (for moment inequalities purpose).

\item[$\bullet$] The second reason is that MDS is not an abstract notion. Indeed, there are a lot of well known
stationary MDS models which are used in practice, such as
ARCH(1) or more general GARCH(1,1) stochastic volatility  models.
\end{itemize}

We compare, in the case of nonparametric regression model with MDS errors, the behaviors of the smoothing bandwidths obtained by minimizing
{\textcolor{black}{either the unknown average squared error, the theoretical mean average squared error,
a Mallows-type criterion adapted to the dependent case and the family of criteria known as generalized cross validation (GCV) extensions of the Mallows' criterion. We prove that these three minimizers and   those based on the  GCV family are first-order equivalent in probability. We give also a normal asymptotic behavior of the
gap between the minimizer of the average square error and that of the  Mallows-type criterion. This is extended to the GCV family.}}
The obtained results generalize those under
independent errors, as in  \cite{Hardele}, to MDS ones.
Finally, we apply our results to a specific case of MDS namely the ARCH(1) processes.

The adaptation to the dependent case from the independent one is not trivial and  needs to establish more theoretical and technical results such as maximal inequalities or limit theorems for quadratic forms of dependent data.
To establish our theoretical results, we make use of some ingredients adapted to our case of dependent observations taken from \cite{bur}, \cite{DLLY},
\cite{Mc} and \cite{Rio}. {\textcolor{black}{Those ingredients are stated in Appendices {\bf{B}} and {\bf{C}} of the supplementary material. Their}} proofs are based, in particular, on {\textcolor{black}{Marcinkiewicz-Zygmund type inequalities}} or maximal moment inequalities for {\textcolor{black}{MDS, weighted sums of MDS or quadratic forms for MDS}} that we establish using Burkholder-type moment inequalities
together with some chaining arguments ({\textcolor{black}{see Lemma C.1.,$\cdots$, Lemma C.4 and Theorem C.1., Corollary C.1. and Proposition C.1. of Appendix {\bf{C}}}}). Recall that chaining is a nice  approach to approximate the supremum, over a non countable set, of stochastic processes
(used in the theory of empirical processes see for instance \cite{AP}, \cite{SL},  or \cite{DP}).  A central limit theorem for triangular arrays of quadratic forms for MDS is also needed for the proofs of our results.
We prove this central limit theorem, {\textcolor{black}{in Appendix {\bf{B}}}}, by checking the technical conditions of \cite{Mc}.

Our paper is organized as follows. In Section \ref{sec2}, we introduce the regression model and the different  criteria for the selection
of the smoothing parameter $h$. In Section \ref{result}, we state our main results. We apply our theoretical results, in Subsection \ref{arch}, to  ARCH(1) processes.
A Monte-carlo simulation study is conducted in Subsection \ref{arch1}.
%The proofs of our results are given in Section \ref{sec4}. Appendices {\bf{A}} and {\bf{B}} of the supplementary material are dedicated to the proofs of
%the main tools needed to establish our main results.
%Appendix {\bf{C}} of the supplementary material gives and proves some ingredients for MDS used throughout the proofs of the main results {\textcolor{black}{(such as Marcinkiewicz-Zygmund type inequalities or maximal bounds for weighted sums of MDS or quadratic form of MDS as stated above)}}.
%
An essential part of the proofs of our results are givne in Section  \ref{sec4}.
A supplementary material \cite{KDSsupp} is provided containing the rest of the proofs required for the results of this paper as well as some technical ingredients (see section Supplement).

\section{Model and notations}\label{sec2}

Let $(\epsilon_i)_{i\geq 0}$ be a stationary sequence of centered random variables with finite second moment.
Let $\sigma^2=\Var(\epsilon_1)$ and $R$ be the correlation matrix of the vector
$(\epsilon_1,\cdots,\epsilon_n)$.
Consider the following regression model, defined  for $i=1,\cdots,n$, by
\begin{equation}\label{model}
Y_i=r(x_i)+\epsilon_i,\,\,\,\,\, x_i=\frac{i}{n},
\end{equation}
where $r$ is an unknown regression function  of class ${\cal C}^2$ and the $x_i$'s are equally spaced fixed design.
We are interested in this paper by the Priestley-Chao estimator of $r$ defined,  for $x\in \BBr$, by
$$
{\hat r}(x)=\sum_{i=1}^n l_i(x) Y_i,\,\,\, {\mbox{with}}\,\,\,\,\, l_i(x)= \frac{1}{nh}K\left(\frac{x-x_i}{h}\right),
$$
where
 $K$ is a compactly supported even kernel with class {\textcolor{black}{${\cal C}^2([-1,1])$}} and $h$ is a positive bandwidth less than $1/2$.
The above curve estimator entails the following smoothing,  in the matrix form,
$$
{\hat r}= LY \,\,\, {\mbox{with}}\,\,\,\,\,{\hat r}= ({\hat r}(x_1),\cdots, {\hat r}(x_n))^t
,\,\,\, Y=(Y_1,\cdots,Y_n)^t$$
and  $L= (l_j(x_i))_{1\leq i, j\leq n}$ is known as the smoothing matrix or the
hat matrix.
Since the estimator ${\hat r}$  depends on some smoothing parameter $h$,
we will need some procedure for choosing $h$. For this, we  recall some known criteria of selecting this parameter $h$.

In order to eliminate the boundary
effects of the compactly supported kernel $K$, we introduce, as was done in the literature (see for instance \cite{GM}), a known function
 supported on a sub-interval of the unit interval.
For this, suppose without loss of generality that $h<\epsilon$ where $\epsilon$ is a fixed positive real number less than $1/2$.
Let $u:=u_{\epsilon}$ be a positive function, of class ${\cal C}^1$ and $[\epsilon, 1-\epsilon]$-compactly supported {\textcolor{black}{satisfying $\int_0^1 u(x){{(r^{\prime\prime}}(x))^{2}}dx\neq 0$}}.
Define the average squared error $$
T_n(h)=\frac{1}{n}\sum_{i=1}^n u(x_i)(\hat r(x_i) -r(x_i))^2= \frac{1}{n}\|U^{1/2}(\hat r -r)\|^2,
$$
where $U$ is the diagonal matrix $U=diag(u(x_1),\cdots,u(x_n))$ and for any vector $v$, $\|v\|^2=v^tv$.
\\
\textcolor{black}{It should be pointed out that in order to overcome the boundary problem mentioned above, one may consider the local linear estimate, as done for instance by
\cite{FG}, \cite{FOV}, \cite{BD} where plug-in asymptotic methods for selecting the smoothing parameter have been considered for some class of correlated errors. However the purpose of the current work is the study of bandwidth selection methods which are based on unbiased (or nearly unbiased) criteria for any fixed sample size.  We believe that the extension of our results to local linear estimate may be carried out but it is beyond the scope of this paper and thus could be treated in a separate possible future work.
}

The following lemma (its proof is given in Appendix {\bf{A.1}} of the supplementary material) evaluates its mean, $\BBe(T_n(h))$, {\textcolor{black}{for finite variance
stationary errors}} $(\epsilon_i)_{i\in \BBn}$.
\begin{lem}\label{ase}
Suppose that $\sum_{k=1}^{\infty} k |\Cov(\epsilon_{0},\epsilon_{k})|<\infty$.
Define,
\begin{eqnarray*}
&& D_n(h)=\frac{h^4}{4} \int_{0}^1 u(x) {{(r^{\prime\prime}}(x))^{2}} dx \left(\int_{-1}^{1} t^2K(t)dt\right)^2 \\
&& +
\frac{1}{nh} (\int_0^1u(x)dx) \int_{-1}^1 K^2(y)dy \left(\sigma^2+ 2\sum_{k=1}^{\infty}\Cov(\epsilon_0,\epsilon_k)\right).
\end{eqnarray*}
Then for any $n\geq 1$ and $h\in ]0,\epsilon[$,
\begin{eqnarray*}
&& \BBe(T_n(h))
= D_n(h) + O(\frac{1}{n})+o(h^4)+ O(\frac{1}{n^2h^4})+ \frac{\gamma(h)}{nh},
\end{eqnarray*}
where $O$ is uniformly on $n$ and $h$,  $\gamma(h)$ depends on $h$ (but not on $n$) and tends to $0$ when $h$ tends to $0$.
\end{lem}
Let {\textcolor{black}{$h_n^*\in \argmin_{h>0}D_n(h)$}}.  Clearly, {\textcolor{black}{since we supposed that $\int_{0}^1 u(x) {{(r^{\prime\prime}}(x))^{2}} dx\neq 0$}},
$$
h_n^*= n^{-1/5}
\left(\frac{(\int_0^1u(x)dx) \int_{-1}^1 K^2(y)dy \left(\sigma^2+ 2\sum_{k=1}^{\infty}
\Cov(\epsilon_0,\epsilon_k)\right)}{ \int_{0}^1 u(x) {(r^{\prime\prime}}(x))^{2}dx (\int_{-1}^{1} t^2K(t)dt)^2}\right)^{1/5}=: cn^{-1/5}.
$$
Let, as in  \cite{hall} and \cite{roce}, $H_n$ be a neighborhood of $h_n^*$, i.e, $H_n=[an^{-1/5}, bn^{-1/5}]$ for some fixed $a<c<b$. Define also,
$$
{\textcolor{black}{h_n \in \argmin_{h\in H_n}\BBe(T_n(h))\,\,\,{\mbox{and}}\,\,  \hat h_n\in  \argmin_{h\in H_n}T_n(h).}}
$$

Of course these three ``optimal'' parameters $h_n$, $h_n^*$ and $\hat h_n$ depend on the unknown function $r$, since the criteria that they respectively minimise,
depend themselves  on the regression function $r$. Many authors agree that, among these ones, $\hat h_n$ should be the target
(see \cite{DG}, page 316).
For this reason, an important literature considered minimizers of  ``good'' estimators
of $T_n(h)$ and studied their asymptotic behavior.

For i.i.d. errors $(\epsilon_i)_{1\leq i\leq n}$ with all finite moments, this question is solved.
A  { \textcolor{black}{reasonably good estimate
of
$\BBe(T_n(h))$
is constructed}}
allowing to define a criterion
{ \textcolor{black}{ (of course the ``goodness'' can be measured up to a multiplicative positive factor or an additive constant)}}
that selects an observable choice for $h$ :
{ \textcolor{black}{  cross-validation is often used or the following simpler criterion is also used}}
\begin{eqnarray}\label{Mallows}
&& \hat C_p:=\hat C_p(h) = \frac{1}{\sum_{i=1}^nu(x_i)}\sum_{i=1}^n u(x_i)(Y_i- \hat r(x_i))^2+ 2\frac{\nu}{n}\hat \sigma_h^2,
\end{eqnarray}
where,
$$
\hat \sigma_h^2 := \frac{1}{\sum_{i=1}^nu(x_i)}\sum_{i=1}^n u(x_i)(Y_i- \hat r(x_i))^2 \,{\mbox{ and }} \, \nu:=n\frac{tr(UL)}{tr(U)}=\frac{1}{h}K(0).
$$
The above notation $\hat C_p$, where $\nu$ is the { \textcolor{black}{``local'' (or more generally ``weighted'')}} ``degrees of freedom'',  is related to the $C_p$-statistics introduced by   \cite{MA} for
variable selection in linear regression models.
{ \textcolor{black}{ {Notice that
$\hat C_p(h)  = \hat \sigma_h^2  \times  \Xi_{\rm S} (t(h))$  with   $\Xi _{\rm S}(t):=1+2t$ and  $t(h) := \nu/n$, that means that $\hat C_p(h)$ coincides with the Shibata criterion, as named by  \cite{Hardele}. Recall that the main result of these authors was that they showed a second-order equivalence
(defined below as a footnote)  of the ``exact'' $C_p$ criterion, which uses the exact $\sigma$ instead of $\hat\sigma_h$ in (2.2), and
any criterion obtained by replacing in this second expression of $\hat C_p$ the  penalization factor function $t \longmapsto\Xi_{\rm S}(t)$ by any function $\Xi_{\rm X} $ satisfying
\begin{eqnarray}\label{penalization}
\Xi_{\rm X} (t)= 1 + 2 t + O(t^2)    {\mbox{ with second derivative }}   \Xi_{\rm X}^{\prime\prime}  {\mbox{ bounded on a neighborhood of }}   0;
\end{eqnarray}
a second example being the popular GCV criterion,
$ {\hat \sigma_h^2} / { \left(  1- t(h) \right)^2}   $, associated with the choice $\Xi_{\rm X} (t):=    \Xi_{\rm GCV} (t)  =  (1-t)^{-2}$.
}}
{\textcolor{black}{ Precisely, letting}} $\hat h$ be a minimizer over $h\in H_n$ of the exact $C_p$ criterion,  \cite{Hardele}
proved, in the context of  i.i.d errors $(\epsilon_i)_{1\leq i\leq n}$ with all finite moments,  that $\hat h, h_n^*, \hat h_n, h_n$
are all equivalent in probability,
that  $\hat h-\hat h_n$, $h_n-\hat h_n$ are also close in distribution
as $n$ tends to infinity,
{\textcolor{black}{ and that minimizing any such criterion $\hat \sigma_h^2 \times \Xi_{\rm X} (t(h))$ also produces a bandwidth which is second-order equivalent to $\hat h$}}.
\footnote{The second order equivalence of the ${\rm C_p}$ and ${\rm GCV}$ selectors means that the asymptotic law of $\hat h-\hat h_n$ is unchanged if $\hat h$ is replaced by the minimizer of ${\rm GCV(h)}$}.

The above criteria can hardly be considered as adapted  to the case of
{ \textcolor{black}{  general}}
dependent errors since they take into account
only the variance $\sigma^2$ of the errors and not  their overall dependence
structure. Several authors  extended Mallows'  criterion to some cases of stationary dependent errors.
 \cite{YW} and  \cite{HG}, among others,
generalized Mallows'  criteria  in  (\ref{Mallows}) (but for other purposes than ours) to stationary dependent errors with known covariance matrix   $\sigma^2R$ of the vector
$(\epsilon_1,\cdots,\epsilon_n)^t$, by
\begin{eqnarray}\label{Mallowsdep}
&& {\rm CL}(h)=n^{-1}\|U^{1/2}(I-L)Y\|^2+ 2\sigma^2 n^{-1} tr(URL),
\end{eqnarray}
which is linked to the average squared error $T_n(h)$ due to the following relation,
\begin{eqnarray*}
&& {\rm CL}(h)= T_n(h)+ \delta_2(h)+ n^{-1}\|U^{1/2}(Y-r)\|^2,
\end{eqnarray*}
where
\begin{eqnarray}\label{delta2}
&& \delta_2(h)= 2n^{-1}(Y-r)^t U(r-\hat r)+ 2\sigma^2 n^{-1} tr(URL).
\end{eqnarray}
Let us consider, according to our purpose, $\hat h_{M}$
to be the minimizer of the dependent version of the Mallows criterion (\ref{Mallowsdep})
$$
{\textcolor{black}{\hat h_{M}\in \argmin_{h\in H_n}{\rm CL}(h).}}
$$

Recall that we are interested in the problem of selecting the parameter $h$ when the errors form a sequence of stationary and dependent random variables.
As we mentioned in the introduction, we consider through all this paper,  the above regression model with stationary MDS
errors (defined in Conditions (C) of Section \ref{result} below).
Since MDS is a sequence of non-correlated and centered random variables, $R$,  which represents the correlation of the errors, is nothing else but
the identity matrix.
 { \textcolor{black}{
Since $R=I$, it seems natural to
 consider again the substitution of the true  $\sigma$ in (\ref{Mallowsdep})  by the same estimate  $\hat \sigma_h^2$ used above in $\hat C_p$
and to
ask whether such a substitution still provides good bandwidth selectors under stationary   MDS
errors.
 Thus, we consider the following minimizers, denoted by the generic $\genhGCV$
\begin{eqnarray}\label{GCV}
{\textcolor{black}{\genhGCV \in \argmin_{h\in H_n}G_{\rm X}(h)}}   \,\,\, {\rm where} \,\,\, G_{\rm X}(h)  := n^{-1}\|U^{1/2}(I-L)Y\|^2  \times \Xi_{\rm X} \left(\frac{tr(UL)}{tr(U)}\right)
\end{eqnarray}
where $\Xi_{\rm X}$, satisfying (\ref{penalization}),  is associated with one of the classical GCV-type criteria.
}}

\section{Main results and applications}\label{result}
The following conditions are required to establish our main results. \\
\\
{\bf{Conditions (C).}} Assume that the errors $(\epsilon_i)_{i\geq 0}$ form a  stationary  MDS with respect to some natural filtration $({\cal F}_i)_{i\geq 1}$, i.e,
for any $i>0$, $\epsilon_i$ is ${\cal F}_i$-measurable and $\BBe(\epsilon_i|{{\cal F}_{i-1}})=0$. Suppose also that
$\BBe(\epsilon_1^{2p})<\infty$ for some $p>8$.
\\
\\
Our first result states that for MDS errors, the  bandwidths $h_n,h_n^*, \hat h_n, \hat h_{M}$  { \textcolor{black}{and $\genhGCV$}} are  first-order equivalent in probability
({\textcolor{black}{in other words, both}} the ${\rm CL}$ { \textcolor{black}{and $G_X$}} criteria enjoy the same ``asymptotic optimality" property).

\begin{pro}\label{pro1} Suppose that Conditions (C) are satisfied.  Then \footnote{Here, and for all the evoked argmin, the results apply to any points of the argmin sets}
 $$\frac{h^*_n}{h_n},
\frac{\hat h_n}{h_n}, \frac{\hat h_{M}}{h_n}  { \textcolor{black}{, \frac{\genhGCV}{h_n}}}$$ all converge  in probability to $1$ as
$n$ tends to infinity.
\end{pro}
\noindent Notice that \cite{hall} gave two theorems for two bandwidth selection methods (precisely a block-bootstrap method
and the classical leave-$k$-out technique{\textcolor{black}{, and the mentioned theorems are respectively their Theorem~2.2 and Theorem~2.3}}) under a rather general dependence assumption on the error sequence,
namely the Rosenblatt mixing condition (see their Section 2.2).
Each of these two theorems is a first-order optimality like Proposition \ref{pro1} above, and it could be applied, in particular, to certain stationary  MDS. However we point out that these two theorems also require that all moments of
the marginal law of the errors are finite. Thus the results of \cite{hall}
cannot be applied  to any ARCH process except the trivial one ($\alpha=0$ in the notation of  Section \ref{arch}). \\
\noindent Our second result  gives, under  {\textcolor{black}{a block-covariance decay}} condition, the rate at which
 $\hat h_n - {\hat h}_{M}$
  { \textcolor{black}{and    $\hat h_n - \genhGCV$}}
  {\textcolor{black}{converge}}
 in distribution to a {\textcolor{black}{common} centered normal law, and furthermore states that the martingale difference dependence does   {\textcolor{black}{\it not}} impact
 this law.

\begin{thm}\label{main} Suppose that Conditions (C) are satisfied.
Moreover, suppose that there exists a positive decreasing function $\Phi$ defined on $\BBr^+$ satisfying
$$
\sum_{s=1}^{\infty}s^4\Phi(s)<\infty,
$$
and for any positive integer $q \leq 6$, $1\leq i_1\leq \cdots\leq i_k<i_{k+1}\leq \cdots\leq i_q\leq n$ such that
$i_{k+1}- i_k \geq \max_{1\leq l\leq q-1}(i_{l+1}-i_l)$,
\begin{eqnarray}\label{cov1}
&& |\Cov(\epsilon_{i_1}\cdots\epsilon_{i_k}, \epsilon_{i_{k+1}}\cdots\epsilon_{i_q})|\leq \Phi(i_{k+1}- i_k),
\end{eqnarray}
{\textcolor{black}{where $\epsilon_{i_1}\cdots\epsilon_{i_k}$ denotes the product $\prod_{\ell=1}^{k}\epsilon_{i_\ell}$ (and likewise for $\epsilon_{i_{k+1}}\cdots\epsilon_{i_q}$)}}.
Then    {\textcolor{black}{both }}
$$
n^{3/10}(\hat h_M - {\hat h}_{n})   { \textcolor{black}{\,\,\, {\rm and} \,\,\,   n^{3/10}(\genhGCV - \hat h_n)  }}
$$
  {\textcolor{black}{converge}}
   in distribution to a centered normal law with variance $\Sigma^2$ given by
\begin{eqnarray*}
\Sigma^2=\frac{4\sigma^{6/5}}{5^2 A^{8/5}B^{2/5}}
\left(\left(\int t^2K(t)dt\right)^2 \int_0^1\hspace{-0.2cm} u^2(x)r^{\prime\prime 2}(x)dx
+
\frac{2 A}{B}
\int_0^1\hspace{-0.2cm} u^2(x)dx\int (K-G)^2(t)dt  \right),
\end{eqnarray*}
where $\sigma^2=\BBe(\epsilon_1^2)$, $G$ is the function defined for any $x\in \BBr$ by $G(x)=-xK'(x)$ and
\begin{eqnarray*}
&& A= \int_0^1 u(x)r^{\prime\prime 2}(x)dx \left(\int t^2K(t)dt\right)^2,\,\,\, B=\int_0^1 u(x)dx  \int K^2(t)dt.
\end{eqnarray*}
\end{thm}
\noindent
\begin{rem}\label{rem1}
The control of the covariance quantity $|\Cov(\epsilon_{i_1}\cdots\epsilon_{i_k}, \epsilon_{i_{k+1}}\cdots\epsilon_{i_q})|$ appearing in (\ref{cov1})
is well known in the literature. It was used, for instance in \cite{DLLY}, in order to obtain Marcinkiewicz-Zygmund
type moments inequalities of an even order of the partial sum $\sum_{i=1}^n\epsilon_i$. If the sequence $(\epsilon_n)_n$ is strongly mixing with
mixing coefficients $(\alpha_s)_{s\in \BBn}$,
then it is proved by \cite{Rio}, see also Lemma 9 in \cite{DLLY}) that,
for $1\leq i_1\leq \cdots\leq i_k<i_{k+1}\leq\cdots\leq i_q\leq n$ such that $s:=i_{k+1}- i_k \geq \max_{1\leq l\leq q-1}(i_{l+1}-i_l)$,
$$
|\Cov(\epsilon_{i_1}\cdots\epsilon_{i_k}, \epsilon_{i_{k+1}}\cdots\epsilon_{i_q})| \leq 4 \int_0^{\alpha_s}Q^q(u)du,
$$
where $Q$ is the quantile function of $|\epsilon_1|$, i.e. the inverse of the tail function $t \longmapsto \BBp(|\epsilon_1|>t)$.
\end{rem}

%\noindent
%\begin{rem}\label{rem2}
%In the i.i.d framework, based on the cross-validation criteria and using local linear estimator, \cite{LR} obtained some different asymptotic results than Theorem \ref{main} The extension  of this theorem to local linear estimator for dependent MDS errors is somewhat interesting  and may be carried out but it is beyond the scope of this paper.
%\end{rem}

\subsection{Application to ARCH(1) processes}\label{arch}
We consider the regression model defined in (\ref{model}) with an ARCH(1) error process $(\epsilon_n)_{n\geq 1}$ defined, for $n\geq 1$, by the following stochastic difference equation,
\begin{equation}\label{earch}
\epsilon_n= \eta_n\sqrt{\sigma^2(1-\alpha)+\alpha \epsilon^2_{n-1}},\,\,0\leq \alpha<1,\,\, \sigma^2>0
\end{equation}
where $(\eta_n)_{n\geq 1}$ is an i.i.d.~ centered sequence distributed as a standard normal law and such that $\eta_n$ is independent of
$(\epsilon_{1},\cdots,\epsilon_{n-1})$.
\begin{pro}\label{proarch} Let $(\epsilon_n)_{n\geq 1}$ be a strictly stationary ARCH(1) process {\textcolor{black}{satisfying}} (\ref{earch}) with $\alpha$ such that
$\alpha^8\prod_{i=1}^8(2i-1)<1$ (this is equivalent to $\alpha< 2025027^{-1/8}  \approx 0.162796$).
Then the conclusions of Proposition \ref{pro1} and Theorem \ref{main} hold.
\end{pro}
\noindent {\bf{Proof of Proposition  \ref{proarch}.}}
 We first recall the following well known properties in the literature (see for instance \cite{engle}, \cite{AML} and the references therein).
\begin{lem} Consider the process $(\epsilon_n)_n$ as defined in (\ref{earch}).  Then
\begin{enumerate}
\item $(\epsilon_n)_n$ is a geometric ergodic homogeneous Markov chain with a unique stationary distribution $\pi$. The stationary distribution $\pi$
is continuous and symmetric.
\item $(\epsilon_n)_n$ is strongly mixing with mixing coefficients $(\alpha_l)_{l>0}$
$$
\alpha_l:=\sup_{A\in \sigma(\epsilon_s,\,s\leq 0),\,B\in \sigma(\epsilon_s,\,s\geq l)}|\Cov(\BBone_A,\BBone_B)|=O(\rho^l),
$$
for some $\rho\in ]0,1[$. {\textcolor{black}{Here $\sigma(\epsilon_s; s \leq 0)$ and $\sigma(\epsilon_s; s \geq l)$ denote, respectively the sigma-fields generated by $\epsilon_s,$ for $s \leq 0$ and $\epsilon_s,$ for $s \geq l$.}}
\item $\BBe(\epsilon_1^{2r})<\infty,$ for $r\in \BBn\setminus{\{0\}}$, if and only if $\alpha^r \prod_{i=1}^r(2i-1)<1$.
\item $\BBp(|\epsilon_1|>x)\sim c x^{-\kappa}$ as $x$ tends to infinity {\textcolor{black}{(in all this paper the notation $a(x)\sim b(x)$ means that
$\lim_{x\rightarrow\infty}\frac{a(x)}{b(x)}=1$)}},
for some $c>0$ and $\kappa$ is given as the unique positive solution to $\alpha^{\kappa/2}\BBe(|\eta_1|^{\kappa})=1.$
\end{enumerate}
\end{lem}
Letting ${\cal F}_i= \sigma(\eta_{1},\cdots,\eta_{i})$, then $\epsilon_i$ is ${\cal F}_i$-measurable and
\begin{eqnarray*}
&& \BBe(\epsilon_i|{\cal F}_{i-1})=\sqrt{\sigma^2(1-\alpha)+\alpha \epsilon^2_{i-1}} \BBe(\eta_i|{\cal F}_{i-1})= 0.
\end{eqnarray*}
The sequence $(\epsilon_n)_n$ is then a martingale-difference. {\textcolor{black}{Moreover, since}} it is strongly mixing with $\alpha_s\leq C \rho^s$, we get, from Remark \ref{rem1}, {\textcolor{black}{the bound}} (\ref{cov1}), for any given $q$, by using
$$
\Phi(s):=4\int_0^{\alpha_s}Q^q(u)du.
$$
{\textcolor{black}{Our task now is to prove that $\sum_{s=1}^{\infty}s^4\Phi(s)<\infty$}}.
We deduce from $\BBp(|\epsilon_1|>x)\sim c x^{-\kappa}$ as $x$ tends to infinity that $Q(u)=O(u^{-1/\kappa})$, ({\textcolor{black}{$q$ satisfies necessarily}} $q<\kappa$, since $\BBe(|\epsilon_1|^q)<\infty$) and
$$
\int_0^{\alpha_s}Q^q(u)du\leq \int_0^{C \rho^s}Q^q(u)du=O(\rho^{s(1-q/\kappa)}).
$$
Consequently, for some positive constant $C$,
$$
\sum_{s=1}^{\infty}s^4\Phi(s)\leq 4\sum_{s=1}^{\infty}s^4\int_0^{C \rho^s}Q^q(u)du\leq C \sum_{s=1}^{\infty}s^4\rho^{s(1-q/\kappa)}<\infty,\,\,\, {\mbox{since}}\,\,\rho\in]0,1[.
$$

\subsection{A Monte-carlo simulation study for a  ``trend plus ARCH(1) process'' }\label{arch1}

 We  report here on  rather extensive   experiments with ARCH(1) noise and a single example of regression function (called the ``deterministic trend'' here)  taken from \cite{WLC}, and we focus on the questions of how accurate is the approximation provided by  Theorem \ref{main}, and whether the restriction $p>8$ in Conditions (C) that we have required  could be relaxed. These questions were also studied, in the first arXiv version of this article, for another trend, namely the well known ``bell shaped'' example much studied since \cite{roce}, which is a lot smoother trend than the one used here. And for the sake of place we only do this experimental study  for the Mallows criterion {\textcolor{black}{(the possible GCV-like  criteria, satisfying (\ref{penalization}), being rather numerous, cf.~\cite{Hardele})}.
 We choose a noise level for which the noise-to-signal ratio is ``moderate'', precisely $0.32^2$. So, the chosen
trend function is
$$r (x) = c_0 + c_1 \left(   \sin (8 x -4) + 2 \exp \left( -256 (x - 0.5)^2\right) \right),$$ where we add the constants $c_0, c_1$ to the definition by \cite{WLC} only so that the range of $r(x)$ is exactly $[0,1]$ when $x \in [0,1]$): a plot of $r$ is inserted  in the left panel of Figure~1.
Each data set is thus the sum of this trend $r$ evaluated at $x_i=i/n, i=1,\cdots,n,$ plus an  ARCH(1) sequence with a ``persistence''  parameter $\alpha$ as defined in the above Subsection \ref{arch}. We consider $6$ settings for  the ARCH(1) noise,  precisely
 $$\alpha \in \{ 0.01,0.162, 0.577, 0.75, 0.9, 0.98 \} {\mbox {, with a common value }} \sigma=0.32.$$
 The first value of $\alpha=0.01$  corresponds nearly to  i.i.d.~normal observation noises (this setting will be referred to as the ``quasi-iid-normal'' case) and the last one generates noise sequences for which a strong serial correlation is always present when the sequence is squared.
  Recall that the moment of order $16$ no longer exists as soon as $\alpha$ is slightly above $0.162$, but the moment  of order
  $4$ still exists for $\alpha < \sqrt{1/3} \approx 0.57735$.

  \smallskip
\begin{center}
\includegraphics[width=13.3cm]{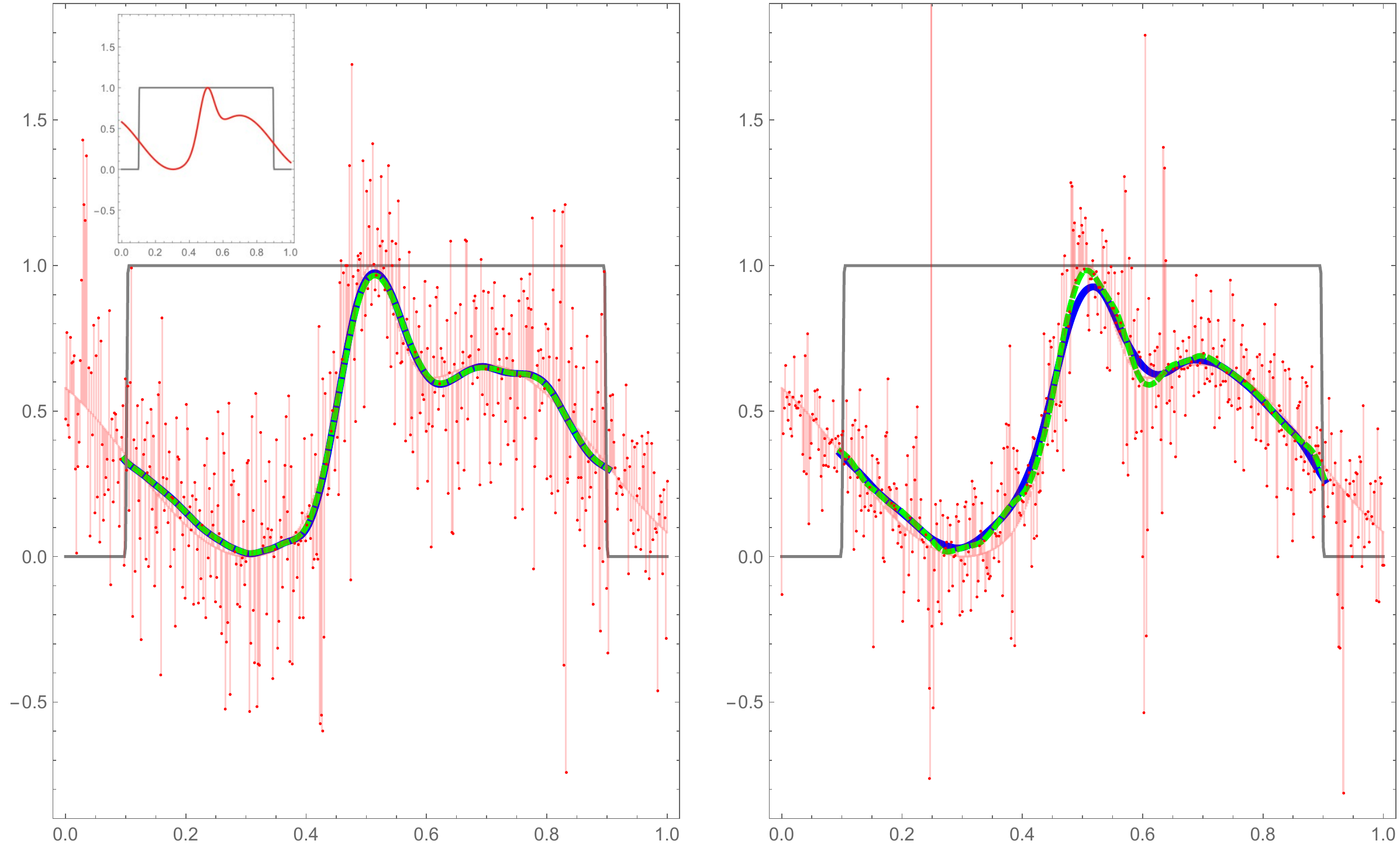}
\end{center}
{\centering \small {Figure 1 : $n=2^{9}$. Each of these  $2$ panels displays one data set $Y$. The underlying  trend $r$ is displayed (in red) in the upper left inset. The  function $u$ (in gray)  excludes about $20\%$ of the points.
The  $2$ panels only differ by $\alpha=0.577$ (left) and $\alpha=0.9$ (right).
For each $Y$ the CL choice and the ``weighted $L_2$-optimal'' $T_n$ choice are plotted  in blue and  dashed green respectively.}}
%\bigskip

\smallskip
The kernel function  used here is the classical biweight
{\textcolor{black}{$K(x) = \frac{15}{16} (1-x^2)^2 1_{[-1,1]}(x)$}}.
%$K(x) = (15/16) (1-x^2)^2 1_{[-1,.1]}(x)$.
As is well known, its precise specification, among possible positive ``bell shaped''  kernels,
has a weak impact on the behavior of bandwidth selection techniques.
As weight function $u$, we used a slightly smoothed version of
$\BBone_{[0.1, 0.9]}$. Since its precise specification also has a weak impact, we only give a plot of the
used  $u$ in Figure~1.
Notice that for the considered data sizes $n$ here, it turns out that it is sufficient to consider  {\textcolor{black}{only bandwidths
 that are  lower than $0.1$,}} as candidate bandwidths.
Then, {\textcolor{black}{since the bandwidth $h$ is one-half the support length of $\frac{1}{h} K(\frac{\cdot}{h}),$}} it can be checked that  computing the sub-vector of the Priestley-Chao estimator $L Y$ whose components are restricted to the $x_i$'s in $[0.1, 0.9]$, can always be done by discrete Fourier transforms.
This remark makes affordable the following simulation study even for quite large $n$.

The data sets size $n$ was chosen in  $\{ 2^9, 2^{12}, 2^{15} \}=  \{ 512, 4096, 32768 \}$.
We generated $1000$ replicated data sets for each of these $3 \times 6$ settings.
For each data set, the minimizer of  $T_n(h)$ and the one of ${\rm CL}(h)$ were numerically computed by a simple grid-search over the domain
 {\textcolor{black}{$[ 0.025, 1]\times 10^{-1}$ (notice the ``no smoothing $h$'' is $n^{-1} K(0) \approx 0.0018$ for $n=512$ and when $h$ comes close to $1/2$ one averages over the entire sample, see \cite{Hardele})}},
the grid-step being chosen fine enough so that the ``granularity'' in the $2000$ computed $h$'s has a very weak impact on the conclusions.

\begin{center}
\includegraphics[width=13.cm]{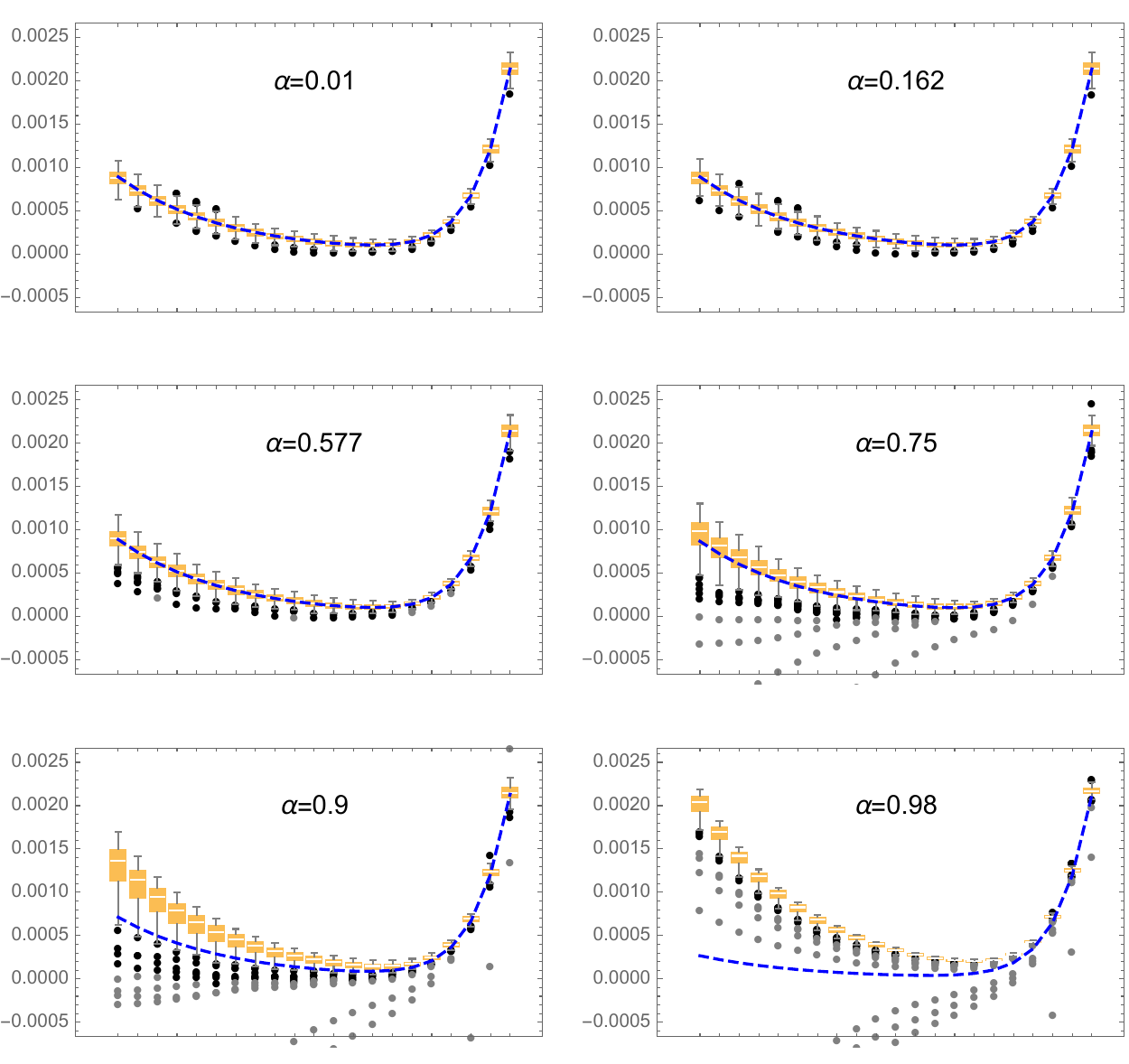}
\end{center}
{\centering \small {Figure 2 : $n=2^{15}$. These $6$ panels only differ by $\alpha$ varying in
 $\{ 0.01,0.162, 0.577, 0.75, 0.9, 0.98 \}$. In each panel, the dashed blue curve is ``empirical MASE'', precisely the average (over the $1000$ replicates)
 of the $T_n(h)$ curves. Each of the 21 boxplots (located at 21 discrete values for $h$ equispaced over {\textcolor{black}{$[ 0.025, 1]\times 10^{-1}$}}) are built from the $1000$ replicates
 of  ${\rm CL}(h) - n^{-1}\|U^{1/2}(Y-r)\|^2$.}}

\bigskip

\noindent
{\bf The ``a.o."~property.}
First, let us analyze the asymptotic optimality (a.o.) result. As is well known, a result like Proposition \ref{pro1}
generally stems from a uniform relative accuracy result which states that ${\rm CL}(h) - n^{-1}\|U^{1/2}(Y-r)\|^2$ uniformly approximates $T_n(h)$
(or its expectation ${\rm MASE}(h)$) with a small (in probability and in $sup$ norm over the domain of candidate $h$'s) error,
``small'' being defined relatively to    ${\rm MASE}(h)$.

We resume in Figure~2 that a uniform relative accuracy is well observed and, above all, this accuracy
in the case $\alpha=0.162$ is of the same order as the accuracy observed in the quasi-iid-normal case
($\alpha=0.01$). Furthermore, an interesting observation is that this accuracy is not deteriorated when $\alpha=0.577$.
However there is clearly a deterioration for larger  $\alpha$, especially  for $\alpha=0.9$ or $0.98$ where in addition to the increased variability, a large bias is observed.
Figure~1 exhibits  such a bias toward oversmoothing for $\alpha=0.9$.

It can be thus  conjectured that, at least for ARCH(1) processes, the restriction $p>8$ of our Conditions~(C) might be weakened to $p>2$. However,
the poor behavior of ${\rm CL}$ (even with quite large $n$) in cases $\alpha=0.75, 0.9$ or $0.98$, leads us to conjecture that
$p>2$ should be considered as a necessary condition for the a.o.~of ${\rm CL}$ or GCV under general  {\textcolor{black}{stationary}} MDS observation errors.

\smallskip

\noindent
{\bf Asymptotic normal distribution.} Now, let us look at the usefulness  of the asymptotic normal approximation stated in Theorem \ref{main}. By inspecting Figure~3,
we clearly see, in the left-bottom panel, that this approximation fits very well for $n=2^{15}$ and  $\alpha=0.577$.

\begin{center}
\includegraphics[width=13.5cm]{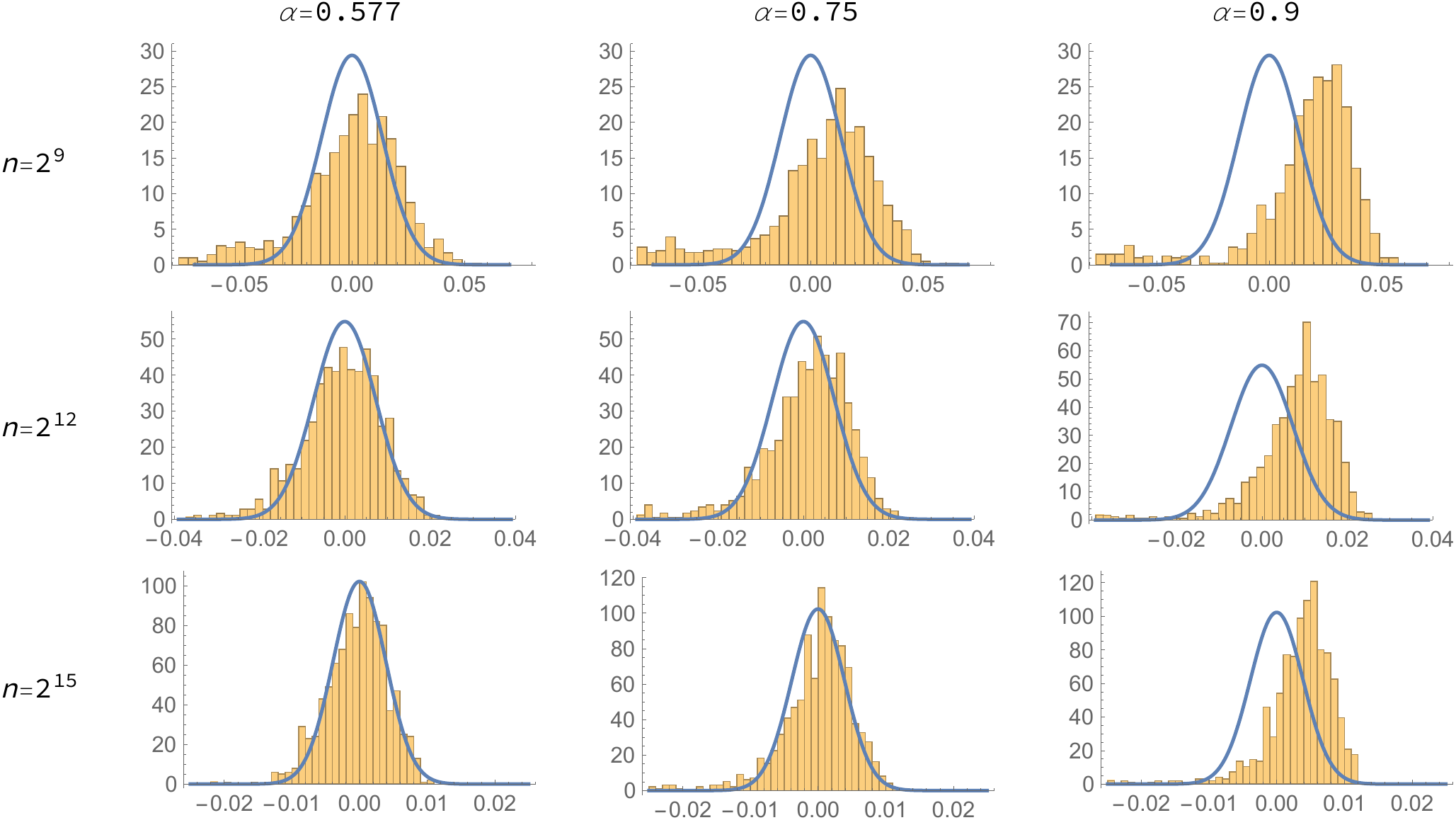}
\end{center}
\bigskip
{\centering \small {Figure 3 : $n=512$ (top), $4096$ (middle) and  $32768$ (bottom). These $9$ panels only differ by $n$ and by $\alpha$
varying in $\{ 0.577, 0.75, 0.9 \}$. In each panel, the displayed histogram is that of the $1000$ replicates of
$\hat h_{M}-\hat h_{n}$; the histograms are normalized so that their integrals are equal to $1$. The superposed blue curve is the normal distribution of $\hat h_{M}-\hat h_{n}$ predicted by the asymptotic theory.
{\textcolor{black}{Notice that, as expected, the range of the abscissae ($h$-differences) decreases by moving from $n=2^{9}$ to  $n=2^{15}$}} .}}
\bigskip

We have also made such a comparison for  $\alpha=0.01$ and  $\alpha=0.162$ (not shown in Figure~3),
and, as expected by Theorem \ref{main} and Section~3.1, the fit is also very good.
For settings with the much smaller $n=512$, the fit is still rather good for  $\alpha=0.577$, but this is no longer true for  $\alpha=0.75$  or  $\alpha=0.9$ (and the fit is even worse for  $\alpha=0.98$, not shown in Figure~3).
{\textcolor{black}{For $n=2^{12}$ we see that the accuracy of the fit, when $\alpha=0.577$, is almost as good as in the case $n=2^{15}$.
Similar conclusions were obtained for the ``bell shaped'' trend mentioned above.}}
It is good news that {\textcolor{black}{the asymptotic approximation given by}} Theorem \ref{main} is thus useful also with $\alpha=0.577$, since this gives support to
the conjecture that Theorem \ref{main} could be extended to an ARCH process under the only existence of the fourth moment of the marginal law.

\section{Proofs}\label{sec4}

\subsection{Main lemmas for the proof of Proposition \ref{pro1}}\label{A.3}
The following two lemmas are very useful for the proof of  Proposition \ref{pro1}.
{\textcolor{black}{Their proofs use the tools (stated in Appendix {\bf{C}}), that control the supremum  of the higher moments of  sums of weighted MDS or quadratic form of weighted MDS, since the evoked quantities $\delta_2(h), T_n(h)$ and their derivatives are expressed in terms of sums of weighted MDS or quadratic form of weighted MDS.}}
\\
 We denote by, $\|\cdot\|_p$  the $p$-norm, i.e, for a random variable $X$, $\|X\|_p= (\BBe(|X|^p))^{1/p}$ and we recall that  $\delta_2(h)$ is defined as in (\ref{delta2}) and $H_n=[an^{-1/5}, bn^{-1/5}]$ for some fixed $a<c<b$.
\begin{lem}\label{P1} It holds, for $p>8$,
\begin{eqnarray}
&& \lim_{n\rightarrow \infty}\|\sup_{h\in H_n}nh | \delta_2(h)|\|_p=0, \label{P1-1} \\
&&
{\textcolor{black}{ \lim_{n\rightarrow \infty}\|\sup_{h\in H_n}nh^2 | \delta_2^{\prime}(h)|\|_p=0, \label{P1-1-firstDeriv} }}\\
&&\lim_{n\rightarrow \infty}\|\sup_{h\in H_n}nh^3 | \delta_2^{\prime\prime}(h)|\|_p=0, \label{P1-3}\\
&& \lim_{M\rightarrow \infty}\limsup_{n\rightarrow \infty}\BBp\left(\sup_{h\in A_{\epsilon}}\sqrt{n} |\delta_2^{\prime\prime}(h)|\geq M\right)=0\label{P1-2},
\end{eqnarray}
where for fixed $\epsilon>0$, $A_{\epsilon}$ is a subset of $H_n$ defined by,
$$
A_{\epsilon}=\{h\in H_n,\,\, \left|\frac{h}{h_n}-1\right|\leq \epsilon\}.
$$
\end{lem}
\begin{lem}\label{P2} It holds, for $p>8$,
\begin{eqnarray}
&& \lim_{n\rightarrow \infty}\|\sup_{h\in H_n}nh |T_n(h)-\BBe(T_n(h))|\|_p=0 \label{P2-1}, \\
&& {\textcolor{black}{  \lim_{n\rightarrow \infty}\|\sup_{h\in H_n}nh^2 |T^{\prime}_n(h)-\BBe(T^{\prime}_n(h))|\|_p=0, \label{P2-1-firstDeriv}}} \\
&& \lim_{n\rightarrow \infty}\|\sup_{h\in H_n}nh^3 |T^{\prime\prime}_n(h)-\BBe(T^{\prime\prime}_n(h))|\|_p=0. \label{P2-2}
\end{eqnarray}
\end{lem}
\subsubsection{Proof of Lemma \ref{P1}}
We have the following decomposition
\begin{eqnarray*}
&& \delta_2(h)=\left(2n^{-1}(Y-r)'U(r-\hat r)+ 2\sigma^2 n^{-1} tr(URL)\right)\\
&&=2\sum_{i=1}^nA_i(h)u(x_i)\epsilon_i +{{2\sum_{i=1}^n\sum_{j=1}^{i-1}(u(x_i)+u(x_j))B_{i,j}(h) \epsilon_i\epsilon_j}}
+ 2\sum_{i=1}^n u(x_i)B_{i,i}(h)\left(\epsilon_i^2-\BBe(\epsilon_i^2)\right),
\end{eqnarray*}
with,
\begin{eqnarray*}
&&{{A_i(h)= \frac{1}{n}\left(r(x_i)- \BBe(\hat r(x_i))\right)=-\frac{1}{n}B(x_i,h)}} \\
&& B_{i,j}(h)= - n^{-1}\frac{1}{nh}K(\frac{x_i-x_j}{h}),\,\,\, B_{i,i}(h)= - n^{-1}\frac{1}{nh}K(0).
\end{eqnarray*}
{\it{Proof of (\ref{P1-1}).}}
We have, for any $h,h'\in H_n$, (using the same calculations yielding to {(\bf{A.11})} of the supplementary material),
\begin{eqnarray*}
&& u(x_i)\sup_{h\in H_n}nh|A_i(h)| =O(n^{-3/5}),\\
&& u(x_i)|nhA_i(h)-nh'A_i(h')|\\
&&\leq hu(x_i)|B(x_i,h)-B(x_i,h')|+ u(x_i)B(x_i,h') |h-h'|\leq cst\, n^{-2/5}|h-h'|,
\end{eqnarray*}
\begin{eqnarray*}
&& nh|B_{i,j}(h)|\leq cst\, n^{-1}\BBone_{|i-j|\leq nh},\\
&& |nhB_{i,j}(h)-nh'B_{i,j}(h')|%\leq \frac{cst}{n}|\frac{1}{h}-\frac{1}{h'}||x_i-x_j|\BBone_{|i-j|\leq n\max(h,h')}
\leq cst\, n^{-4/5}|h-h'|\BBone_{|i-j|\leq n\max(h,h')},
\end{eqnarray*}
and $nh B_{i,i}(h) =\frac{K(0)}{n}$.
Using Lemmas {\bf{C.1}}, {\bf{C.3}} of the supplementary material and the fact that
$$
\frac{1}{n}\sum_{i=1}^nu(x_i)(\epsilon_i^2-\BBe(\epsilon_i^2)) \rightarrow\,0,\, {\mbox{in probability as}}\,\, n\rightarrow \infty,
$$
we get
$$
\lim_{n\rightarrow \infty} \sup_{h\in H_n} nh|\delta_2(h)|=0, \, {\mbox{in probability as}}\,\,  n\rightarrow \infty.
$$
{\it{Proofs of (\ref{P1-1-firstDeriv}) and (\ref{P1-3}).}} {\textcolor{black}{We only discuss the proof (\ref{P1-3}) since that of (\ref{P1-1-firstDeriv}) is similar.}}
For any $h\in ]0,\epsilon[$ and any $n\geq1$,
\begin{eqnarray*}
\delta^{\prime\prime}_2(h) = 2 \sum_{i=1}^nc_i(h)u(x_i)\epsilon_i + 2\sum_{i=1}^n\sum_{j=1}^{i-1} (u(x_i)+u(x_j))c_{i,j}(h) \epsilon_i\epsilon_j
+
2\sum_{i=1}^n c_{i,i}(h) u(x_i)(\epsilon_i^2-\BBe(\epsilon_i^2)),
\end{eqnarray*}
where, letting $K_1=K-G$ and $G_1(u)= -uK_1'(u)$,
\begin{eqnarray*}
&& c_i(h) = -\frac{1}{n}\frac{\partial^2}{\partial h^2} \BBe(\hat r(x_i))= -\frac{1}{n}\frac{\partial^2}{\partial h^2} B(x_i,h), \\
&& c_{i,j}(h) = -\frac{2}{n^2h^3}K_1(\frac{x_i-x_j}{h})+ \frac{1}{n^2h^3} G_1(\frac{x_i-x_j}{h}),
\,\,\, {\rm and} \,\,\,  c_{i,i}(h)=\frac{-2}{n^2h^3}K(0).
\end{eqnarray*}
Now,
\begin{eqnarray*}
&&  nh^3\sum_{i=1}^n c_{i,i}(h) u(x_i)(\epsilon_i^2-\BBe(\epsilon_i^2))=  {cst } \frac{1}{n}\sum_{i=1}^n  u(x_i)(\epsilon_i^2-\BBe(\epsilon_i^2))
\end{eqnarray*}
which converges in probability to $0$ by an analogous to Lemma \ref{lem6}.\\
\\
We also have, for any $h,h'\in H_n$, (see Lemma {\bf{A.3}} of the supplementary material) the following bounds,
\begin{eqnarray*}
&& u(x_i)|nh^3c_{i}(h)|\leq cst\,\, n^{-3/5}, \\
&& u(x_i)|nh^3c_{i}(h)- nh^{'3}c_{i,n}(h')|\leq cst n^{-2/5}|h-h'|, \\
&& u(x_i)|nh^3c_{i,j}(h)|\leq cst\, \frac{1}{n}\BBone_{|i-j|\leq nh},\\
&& u(x_i)|nh^3c_{i,j}(h)-nh^{'3}c_{i,j}(h')| \leq   cst n^{-4/5}|h-h'|\BBone_{|i-j|\leq n\max(h,h')}.
\end{eqnarray*}
All the requirements of Lemmas {\bf{C.1}} and {\bf{C.3}} of the supplementary material are satisfied. We deduce  that,
$$
\lim_{n\rightarrow \infty} \sup_{h\in H_n} nh^3|\delta_2^{\prime\prime}(h)|=0,\,\,\, {\mbox{in probability}}.
$$
{\it{Proof of (\ref{P1-2}).}} We have,
\begin{eqnarray*}
&&
\sqrt{n}\delta^{\prime\prime}_2(h) = 2 \sum_{i=1}^ne_i(h)u(x_i)\epsilon_i + 2\sum_{i=1}^n\sum_{j=1}^{i-1} (u(x_i)+u(x_j))e_{i,j}(h) \epsilon_i\epsilon_j\\
&& +
2\sum_{i=1}^n e_{i,i}(h) u(x_i)(\epsilon_i^2-\BBe(\epsilon_i^2)),
\end{eqnarray*}
where,
\begin{eqnarray*}
&& e_i(h) = -\frac{1}{\sqrt{n}}\frac{\partial^2}{\partial h^2} \BBe(\hat r(x_i))= -\frac{1}{\sqrt{n}}\frac{\partial^2}{\partial h^2} B(x_i,h), \\
&& e_{i,j}(h) = -\frac{2}{n\sqrt{n}h^3}K_1(\frac{x_i-x_j}{h})+ \frac{1}{n\sqrt{n}h^3} G_1(\frac{x_i-x_j}{h}),  \\
&& e_{i,i}(h)=\frac{-2}{n\sqrt{n}h^3}K(0).
\end{eqnarray*}
We have, since $\sum_{i=1}^{\infty}|\Cov(\epsilon_1^2, \epsilon_i^2)|<\infty$,
\begin{eqnarray*}
&& \sup_{h\in H_n}\left|\sum_{i=1}^n e_{i,i}(h) u(x_i)(\epsilon_i^2-\BBe(\epsilon_i^2))\right|\leq \frac{cst}{n^{9/10}}\left|\sum_{i=1}^n  u(x_i)(\epsilon_i^2-\BBe(\epsilon_i^2))\right|\\
&& \left\|\sup_{h\in H_n}\left|\sum_{i=1}^n e_{i,i}(h) u(x_i)(\epsilon_i^2-\BBe(\epsilon_i^2))\right|\,\right\|_2^2\leq
\frac{cst}{n^{18/10}}\left\|\sum_{i=1}^n  u(x_i)(\epsilon_i^2-\BBe(\epsilon_i^2))\right\|_2^2\\
&& \leq cst\,\, \frac{n}{n^{18/10}}.
\end{eqnarray*}
Hence,
\begin{equation}\label{P1-3-1}
\lim_{n\rightarrow \infty}\left\|\sup_{h\in H_n}\left|\sum_{i=1}^n e_{i,i}(h) u(x_i)(\epsilon_i^2-\BBe(\epsilon_i^2))\right|\,\right\|_2=0.
\end{equation}
Let, for $h\in H_n$, $h'\in H_n$, ${\tilde e}_i(h)= e_i(h)-e_i(h_n)$, then
\begin{eqnarray*}
&& u(x_i)|{\tilde e}_i(h)|\leq  cst\frac{1}{\sqrt{n}} |h-h_n| \leq cst\, n^{-7/{10}},\\
&& u(x_i)|{\tilde e}_i(h)-{\tilde e}_i(h')|\leq cst\frac{1}{\sqrt{n}} |h-h'|.
\end{eqnarray*}
Hence,
\begin{eqnarray*}
&& \sup_{h\in H_n}\left|\sum_{i=1}^ne_i(h)u(x_i)\epsilon_i\right|\leq
\sup_{h\in H_n}\left|\sum_{i=1}^n(e_i(h)-e_i(h_n))u(x_i)\epsilon_i\right|+ \left|\sum_{i=1}^ne_i(h_n)u(x_i)\epsilon_i\right|.
\end{eqnarray*}
It follows that
\begin{eqnarray*}
&& \left\|\sup_{h\in H_n}\left|\sum_{i=1}^ne_i(h)u(x_i)\epsilon_i\right|\right\|_p\\ && \leq
\left\|\sup_{h\in H_n}\,\left|\sum_{i=1}^n(e_i(h)-e_i(h_n))u(x_i)\epsilon_i\right|\right\|_p+
 \left\|\sum_{i=1}^ne_i(h_n)u(x_i)\epsilon_i\right\|_p.
\end{eqnarray*}
Applying Lemmas {\bf{C.1}} and Corollary {\bf{C.1}}  of  Appendix {\bf{C}} of the supplementary material, we deduce that,
\begin{equation}\label{P1-3-2}
\limsup_{n\rightarrow \infty} \left\|\sup_{h\in H_n}\left|\sum_{i=1}^ne_i(h)u(x_i)\epsilon_i\right|\right\|_p<\infty.
\end{equation}
Let, for $h\in A_{\epsilon}=\{h\in H_n,\,\, |\frac{h}{h_n}-1|\leq \epsilon\}$,   $H_{i,j}(h)= \frac{1}{h^3}(G_1-2K_1)(\frac{x_i-x_j}{h})$
\begin{eqnarray*}
&& {\tilde e}_{i,j}(h)= e_{i,j}(h)-e_{i,j}(h_n)= \frac{1}{n^{3/2}}(H_{i,j}(h)-H_{i,j}(h_n)).
\end{eqnarray*}
Since, $|\frac{\partial H_{i,j}}{\partial h}(h)|\leq cst\, h^{-4} $, then for any $h,h'\in H_n$
$$
|H_{i,j}(h)- H_{i,j}(h')|\leq  n^{4/5} |h-h'|\BBone_{|i-j|\leq n\max(h,h')}
$$
and, for any $h\in A_{\epsilon}$,
\begin{eqnarray*}
&& |{\tilde e}_{i,j}(h)|\leq cst\,  \frac{\epsilon}{n^{9/10}}|\BBone_{|i-j|\leq n\max(h,h')},\\
&& |{\tilde e}_{i,j}(h)-{\tilde e}_{i,j}(h')|\leq cst\,\frac{1}{{n^{7/10}}} |h-h'|\BBone_{|i-j|\leq n\max(h,h')}.
\end{eqnarray*}
We have,
\begin{eqnarray*}
&& \sup_{h\in A_{\epsilon}}\left|\sum_{i=1}^n\sum_{j=1}^{i-1} (u(x_i)+u(x_j))e_{i,j}(h) \epsilon_i\epsilon_j\right|\\
&& \leq \sup_{h\in A_{\epsilon}}\left|\sum_{i=1}^n\sum_{j=1}^{i-1} (u(x_i)+u(x_j))(e_{i,j}(h)-e_{i,j}(h_n)) \epsilon_i\epsilon_j\right|\\
&& +
\left|\sum_{i=1}^n\sum_{j=1}^{i-1} (u(x_i)+u(x_j))e_{i,j}(h_n) \epsilon_i\epsilon_j\right|.
\end{eqnarray*}
Arguing as in Lemma {\bf{C.4}} of the supplementary material, we have
$$
\limsup_{n\rightarrow \infty}
\left\|\sup_{h\in A_{\epsilon}}\left|\sum_{i=1}^n\sum_{j=1}^{i-1} (u(x_i)+u(x_j))(e_{i,j}(h)-e_{i,j}(h_n)) \epsilon_i\epsilon_j\right|\,\right\|_p<\infty,
$$
and by Proposition {\bf{C.1}} of the supplementary material,
$$
\limsup_{n\rightarrow \infty} \left\|\sum_{i=1}^n\sum_{j=1}^{i-1} (u(x_i)+u(x_j))e_{i,j}(h_n) \epsilon_i\epsilon_j\right\|_p<\infty.
$$
Consequently,
\begin{equation}\label{P1-3-3}
\lim_{n\rightarrow \infty} \left\|\sup_{h\in A_{\epsilon}}\left|\sum_{i=1}^n\sum_{j=1}^{i-1} (u(x_i)+u(x_j))e_{i,j}(h) \epsilon_i\epsilon_j\right|\right\|_p<\infty.
\end{equation}
The limit (\ref{P1-2}) is proved by collecting  (\ref{P1-3-1}), (\ref{P1-3-2}) and  (\ref{P1-3-3}).
\subsubsection{Proof of Lemma \ref{P2}}
We can write the following decomposition,
\begin{eqnarray}\label{e5}
&& T_n(h)-\BBe(T_n(h))=
\frac{1}{n}\sum_{i=1}^nu(x_i)\left[\left(\hat r(x_i)-\BBe(\hat r(x_i))\right)^2 - \BBe\left[\left(\hat r(x_i)-\BBe(\hat r(x_i))\right)^2\right] \right] {\nonumber}\\
&& + \frac{2}{n}\sum_{i=1}^nu(x_i)\left(\hat r(x_i)-\BBe(\hat r(x_i))\right)\left(\BBe(\hat r(x_i))- r(x_i)\right){\nonumber}\\
&&= \sum_{j=1}^n C_{j,n}(h)\epsilon_j+ \sum_{j=1}^n\sum_{l=1}^{j-1}B_{j,l}(h)\epsilon_j\epsilon_l+ \sum_{j=1}^nD_{j,n}(h)(\epsilon_j^2-\BBe(\epsilon_j^2)),
\end{eqnarray}
where,
\begin{eqnarray*}
&& C_{j,n}(h)= \frac{2}{n^2h}\sum_{i=1}^nu(x_i)K(\frac{x_i-x_j}{h})B(x_i,h),\,\,\, B(x_i,h)=\BBe(\hat r(x_i))- r(x_i),\\
&& B_{j,l}(h)= \frac{2}{n^3h^2}\sum_{i=1}^nu(x_i)K(\frac{x_i-x_j}{h})K(\frac{x_i-x_l}{h}), \\
&& D_{j,n}(h)= \frac{1}{n^3h^2}\sum_{i=1}^nu(x_i)K^2(\frac{x_i-x_j}{h}).
\end{eqnarray*}
{\it{Proof of (\ref{P2-1}).}}
Let
\begin{eqnarray*}
&& c_{j,n}(h)=nhC_{j,n}(h)= \frac{2}{n}\sum_{i=1}^nu(x_i)K(\frac{x_i-x_j}{h})B(x_i,h)
\end{eqnarray*}
with $B(x_i,h)=\BBe(\hat r(x_i))-r(x_i)= \frac{1}{nh}\sum_{l=1}^nK\left(\frac{x_i-x_l}{h}\right)r(x_l)-r(x_i)$. We get, for $h,h'\in H_n$,
\begin{eqnarray*}
&& c_{j,n}(h)-c_{j,n}(h')=\frac{2}{n}\sum_{i=1}^nu(x_i) \left(K(\frac{x_i-x_j}{h})B(x_i,h)-K(\frac{x_i-x_j}{h'})B(x_i,h')\right).
\end{eqnarray*}
Now, since $K$ is a Lipschitz function,
\begin{eqnarray*}
&& \left|K(\frac{x_i-x_j}{h})B(x_i,h)-K(\frac{x_i-x_j}{h'})B(x_i,h')\right|\\
&& \leq cst\,\left( |B(x_i,h)-B(x_i,h')| + |x_i-x_j|\sup_{h\in H_n}|B(x_i,h)| \frac{|h-h'|}{hh'}\right)\BBone_{|x_i-x_j|\leq \max(h,h')}.
\end{eqnarray*}
We have, for any $h,h'\in H_n$,
\begin{eqnarray*}
&& u(x_i)|B(x_i,h)-B(x_i,h')| \leq cst\,\, n^{-1/5}|h-h'|,
\end{eqnarray*}
and by the proof of Lemma {\bf{A.1}} in Appendix {\bf{A.1}} of the supplementary material,
$
|B(x_i,h)|\leq cst\, h^2.
$
\noindent Hence, for $h,h'\in H_n$,
{{\begin{eqnarray*}
&& |c_{j,n}(h)-c_{j,n}(h')|\leq cst\, n^{-2/5} |h-h'|.
\end{eqnarray*}}
Since $K$ is compactly supported, we have,
{{\begin{eqnarray*}
&& \sup_{h\in H_n}|c_{j,n}(h)|\leq cst \max_{1\leq i\leq n}\sup_{h\in H_n} (h|B(x_i,h)|)=O(n^{-\frac{3}{5}}).
\end{eqnarray*}}
Consequently, we obtain using Lemma {\bf{C.1}} of the supplementary material,
\begin{equation}\label{e7}
\lim_{n\rightarrow \infty} \left\| \sup_{h\in H_n}nh\left|\sum_{j=1}^n C_{j,n}(h)\epsilon_j\right|\right\|_2=0.
\end{equation}
Now, let
\begin{eqnarray*}
&& d_{j,n}(h)= nhD_{j,n}(h)= \frac{1}{n^2h}\sum_{i=1}^nu(x_i)K^2(\frac{x_i-x_j}{h}).
\end{eqnarray*}
We have, $|d_{j,n}(h)|\leq \frac{cst}{n}$ and  $|d_{j,n}(h)-d_{j,n}(h')|\leq  n^{-4/5}|h-h'|.$
Then Lemma {\bf{C.2}} of the supplementary material gives,
\begin{equation}\label{e9}
\lim_{n\rightarrow \infty} \left\|\sup_{h\in H_n} nh\left|\sum_{j=1}^n D_{j,n}(h)(\epsilon_j^2-\BBe(\epsilon_j^2))\right|\right\|_p=0.
\end{equation}
Now, let
\begin{eqnarray*}
&& b_{j,l}(h)= nhB_{j,l}(h)= \frac{1}{n^2h}\sum_{i=1}^nu(x_i)K(\frac{x_i-x_j}{h})K(\frac{x_i-x_l}{h}),
\end{eqnarray*}
we have,
\begin{eqnarray*}
&&b^2_{j,l}(h)=\frac{1}{n^4h^2}\left(\sum_{i=1}^nu(x_i)K(\frac{x_i-x_j}{h})K(\frac{x_i-x_l}{h})\right)^2 \leq \frac{{\textcolor{black}{cst}}}{n^2} \BBone_{|j-l|\leq 2nh}.
\end{eqnarray*}
Our purpose now is to control, for $h,h'\in H_n$, the increment  $|b_{j,l}(h)-b_{j,l}(h')|$.
We have,
\begin{eqnarray*}
&& |b_{j,l}(h)-b_{j,l}(h')|\leq \frac{cst\,\max(h,h')}{nhh'}|h-h'|\leq cst\,\, n^{-4/5}|h-h'|\BBone_{|j-l|\leq 2nh}.
\end{eqnarray*}
Then by Lemma {\bf{C.3}} of the supplementary material, we obtain
\begin{equation}\label{e10}
\lim_{n\rightarrow \infty} \left\|\sup_{h\in H_n} nh\left|\sum_{j=1}^n \sum_{l=1}^{j-1}B_{j,l}(h)\epsilon_j\epsilon_l\right|\right\|_p=0.
\end{equation}
Collecting (\ref{e5}), (\ref{e7}),  (\ref{e9}) and (\ref{e10}), we finally deduce (\ref{P2-1}).
\\
{\it{Proofs of (\ref{P2-1-firstDeriv}) and (\ref{P2-2}).}} {\textcolor{black}{Let us note that the function $h \longmapsto \BBe(T_n(h))$ is twice differentiable with continuous second derivative and that  for $i\in \{1,2\}$
$\frac{\partial ^i}{\partial h^i} \BBe(T_n(h))= \BBe(T^{(i)}_n(h))$. We only discuss the proof (\ref{P2-2}) since that of (\ref{P2-1-firstDeriv}) is similar.}} %
Taking the second derivative over $h$ in (\ref{e5}), we have
$$ T^{\prime\prime}_n(h)-\BBe(T^{\prime\prime}_n(h))
= \sum_{j=1}^n C^{\prime\prime}_{j,n}(h)\epsilon_j+ \sum_{j=1}^n\sum_{l=1}^{j-1}B^{\prime\prime}_{j,l}(h)\epsilon_j\epsilon_l+ \sum_{j=1}^nD^{\prime\prime}_{j,n}(h)(\epsilon_j^2-\BBe(\epsilon_j^2)),
$$
where, letting $B(x_i,h)=\BBe(\hat r(x_i))- r(x_i)$,
\begin{eqnarray*}
&& C^{\prime\prime}_{j,n}(h)
= \frac{1}{n^2h}\sum_{i=1}^nu(x_i) \left(\frac{B(x_i,h)}{h^2}G_1(\frac{x_i-x_j}{h})+\frac{B'(x_i,h)}{h}G_2(\frac{x_i-x_j}{h}) \right. \\ && \left. + B^{\prime\prime}(x_i,h)G_3(\frac{x_i-x_j}{h})\right),\\
&& B^{\prime\prime}_{j,l}(h)= \frac{1}{n^3h^4}\sum_{i=1}^nu(x_i)F_1(\frac{x_i-x_j}{h})F_2(\frac{x_i-x_l}{h}) \\
&& D^{\prime\prime}_{j,n}(h)= \frac{1}{n^3h^4}\sum_{i=1}^nu(x_i)F(\frac{x_i-x_j}{h}),
\end{eqnarray*}
where $F_1$, $F_2$, $F$, $G_1$, $G_2$ and $G_3$ are bounded functions of class $C^1$, $[-1,1]$-compactly supported.
The proof of (\ref{P2-2}) is analogous to (\ref{P1-2}) and (\ref{P2-1}).

\subsection{Proof of Proposition \ref{pro1}}\label{proof1}
We have from Lemma \ref{ase},
\begin{equation}\label{lim1}
\lim_{n\rightarrow \infty}\sup_{h\in H_n} \left|\frac{\BBe(T_n(h))}{D_n(h)}-1\right|=0.
\end{equation}
From this, we claim that $\frac{h_n}{h^*_n}$ converges  to $1$, as $n$ tends to infinity. In fact, by the definition of $h_n$, it holds
$
\BBe(T_n(h_n)) \leq \BBe(T_n(h_n^*)).
$
Hence,
$$
D_n(h_n)\frac{\BBe(T_n(h_n))}{D_n(h_n)} \leq \frac{\BBe(T_n(h_n^*))}{D_n(h_n^*)}D_n(h_n^*),
$$
so  by (\ref{lim1}) and the definition of $h_n^*$, we deduce that, for a fixed $\epsilon>0$ there exists $n_0$ such that for any $n\geq n_0$,
$$
(1-\epsilon)D_n(h_n)\leq (1+\epsilon)D_n(h^*_n) \leq (1+\epsilon)D_n(h_n),
$$
so that $\lim_{n\rightarrow \infty} \frac{D_n(h^*_n)}{D_n(h_n)}=1$, which ensures that $\lim_{n\rightarrow \infty}\frac{h_n^*}{h_n}=1$, in fact (supposing without loss
of generality that $D_n(h_n^*)-D_n(h_n)\neq 0$),
\begin{eqnarray*}
&& \frac{h_n^*-h_n}{h_n}= \frac{h_n^*-h_n}{D_n(h_n^*)-D_n(h_n)}\frac{D_n(h_n^*)-D_n(h_n)}{h_n}\\
&&=  \frac{h_n^*-h_n}{D_n(h_n^*)-D_n(h_n)}\frac{D_n(h_n)}{h_n}\left(\frac{D_n(h^*_n)}{D_n(h_n)}-1\right)
= \frac{D_n(h_n)}{h_nD_n'(h^*)}\left(\frac{D_n(h^*_n)}{D_n(h_n)}-1\right),
\end{eqnarray*}

where $h^*$ is between $h_n$ and $h_n^*$ which are all in $H_n$, consequently $\limsup_{n\rightarrow \infty}|\frac{D_n(h_n)}{h_nD_n'(h^*)}|<\infty$ and then the behavior of
$ \frac{h_n^*-h_n}{h_n}$ is deduced from the fact that  $\frac{D_n(h^*_n)}{D_n(h_n)}-1$  tends to $0$ as $n$ tends to infinity.

In order to complete the proof of Proposition \ref{pro1}, we only need  to prove that both
$\frac{D_n(\hat h_n)}{D_n(h_n)}$  and $\frac{D_n(\hat h_{M})}{D_n(h_n)}$
converge  in probability to $1$ as $n$ tends to infinity (recall that both $\hat h_{M}$ and $\hat h_n$ belong to $H_n$). We refer the reader to \cite{roce}
for similar arguments.
For this, we have
to prove an analogous to the limit (\ref{lim1}),
\begin{equation}\label{lim2}
\lim_{n\rightarrow \infty}\left\|\sup_{h\in H_n} \left|\frac{T_n(h)}{D_n(h)}-1\right|\,\right\|_p = 0,\,\,\, {\mbox{for some}}\,\,p>8
\end{equation}
which gives, from the same previous arguments, that, for any $\epsilon>0$,
$$
\lim_{n\rightarrow \infty}\BBp\left((1-\epsilon)\leq (1+\epsilon)\frac{D_n(h^*_n)}{D_n(\hat h_n)} \leq (1+\epsilon) \right)=1.
$$
Since $\inf_{h\in H_n}nhD_n(h)>0$  and by (\ref{lim1}), the limit (\ref{lim2}) is proved as soon as,
$$
\lim_{n\rightarrow \infty}  \left\|\sup_{h\in H_n}nh \left|T_n(h) -\BBe(T_n(h))\right|\right\|_p=0,
$$
which immediately follows from Lemma \ref{P2} (more precisely (\ref{P2-1}) of Subsection \ref{A.3}).
Our purpose now is to prove that $\frac{\hat h_{M}}{h_n}$ converges in probability to $1$ as $n$ tends to infinity.
Recall that
$
{\rm CL}(h)= T_n(h)+ \delta_2(h)+ n^{-1}\|U^{1/2}(Y-r)\|^2,
$
where $\delta_2(h)=2n^{-1}(Y-r)'U(r-\hat r)+ 2\sigma^2 n^{-1} tr(URL)$.
We have, using Lemma \ref{P1} of Subsection \ref{A.3},
$$
\lim_{n\rightarrow \infty} \|\sup_{h\in H_n} nh|\delta_2(h)|\,\|_p=0,
$$
or equivalently, since $\inf_{h\in H_n}nhD_n(h)>0$,
$$
\lim_{n\rightarrow \infty} \left\|\sup_{h\in H_n} \left|\frac{\delta_2(h)}{D_n(h)}\right|\,\right\|_p =0.
$$
This last limit, together with (\ref{lim2}), give
\begin{equation}\label{limlim}
\lim_{n\rightarrow \infty} \left\|\sup_{h\in H_n} \left|\frac{T_n(h)+ \delta_2(h)}{D_n(h)}-1\right|\,\right\|_p =0.
\end{equation}
Now, we have, since $n^{-1}\|U^{1/2}(Y-r)\|^2$ doesn't depend on $h$,  $${\textcolor{black}{\hat h_{M}\in \argmin_{h\in H_n}\left(T_n(h)+ \delta_2(h)\right),}}$$
so that using (\ref{limlim}) and the same previous arguments, we prove that
$$
\frac{D_n(\hat h_{M})}{D_n(h_n)}\rightarrow 1, \,\, {\mbox{in probability as}}\,\,\, n\rightarrow \infty.
$$

 {\textcolor{black}{It remains to prove the statement concerning the fourth ratio ${\genhGCV}/{h_n}$.
 A way to do it is to appeal to the a.o. of ${\rm CL}$ that we have proved above and to show that the difference $G_X(h) -   {\rm CL}(h)$ is uniformly negligible as compared to  $\BBe(T_n(h))$.
 Using that $tr(U R L) = \sum_{i=1}^n u(x_i) t(h)$ where $t(h) = n^{-1} h^{-1} K(0)$, we have from the definition (\ref{GCV}) of $G_X(h)$ and the property $ \Xi_{\rm X} (t) =1+2t + O(t^2)$
 $$  \frac{n}{\sum_{i=1}^nu(x_i)} \left(G_X(h) -   {\rm CL}(h) \right)$$
$$=    \frac{1}{\sum_{i=1}^nu(x_i)} \|U^{1/2}(I-L)Y\|^2   \times  \Xi_{\rm X} (t(h))  -
\left(   \frac{1}{\sum_{i=1}^nu(x_i)}  \|U^{1/2}(I-L)Y\|^2 + 2 \sigma^2  t(h)    \right) $$
$$=    \hat \sigma_h^2  \times  \left(    \Xi_{\rm X} (t(h)) -1 \right)  -
 2 \sigma^2  t(h)
  =  \left(  \hat \sigma_h^2  - \sigma^2 \right)  ( 2 t(h))   + \hat \sigma_h^2  O(t(h))^2.$$
Clearly this last quantity is $o_P(  n^{-1} h^{-1} )$ uniformly over a domain of $h$ where $n h \rightarrow \infty$ as soon as $\hat \sigma_h^2$ converges toward $\sigma^2$ uniformly over this domain.
To see that this is true for the domain $H_n$ under our assumptions, let us decompose
\begin{eqnarray*}
&&\frac{\sum_{i=1}^nu(x_i)}{n} \hat \sigma_h^2 \\ &&=
\BBe T_n(h) + \left(T_n(h)  -  \BBe T_n(h)  \right)+ \delta_2(h)+ \frac{1}{n}\|U^{1/2}(Y-r)\|^2 -  2\sigma^2 \frac{\sum_{i=1}^nu(x_i)}{n}  t(h),
 \end{eqnarray*}
 (directly obtained by combining (\ref{Mallowsdep}) and (\ref{delta2})).
Now the required convergence is a consequence of Lemma~\ref{ase}, the  limit (\ref{P2-1}) in Lemma \ref{P2},  the  limit (\ref{P1-1}) in Lemma \ref{P1} and the fact that $ \frac{1}{\sum_{i=1}^nu(x_i)}  \|U^{1/2}(Y-r)\|^2 \to \sigma^2$ in probability under Conditions (C).
 }}
The proof of Proposition \ref{pro1} is completed.
\subsection{Proof of Theorem \ref{main}}\label{proof2}
The following   lemma is crucial for the proof of Theorem \ref{main}. It gives conditions under which
$v_n(\hat h_{M}- \hat h_{n})$  {\textcolor{black}{and $v_n(\hat h_{M}- \hat h_{n})$ converge}}}}   to a normal law with some rate $v_n$.
 {\textcolor{black}{Its proof
 is given in Section {\bf {A.3}} of the supplementary material.}}
\begin{lem}\label{lemtec} If, as $n$ tends to infinity and for some positive rate $a_n$,
(recall that $\delta_2(h)$ is defined in Equation (\ref{delta2})),
\begin{enumerate}
\item $a_n\delta_2'({h}_n)$ converges to a centered normal law with variance $V$,
\item $a_n(\delta_2'({h}_n)- \delta_2'({\hat h_n}))$ converges in probability to $0$,
\item $\frac{{\rm CL}''(h^*)}{\BBe(T_n^{\prime\prime}(h_n))}$ tends in probability to $1$ for any $h^*$ between $\hat h_n$ and $\hat h_{M}$,
\end{enumerate}
then
$$
a_n \BBe(T_n^{\prime\prime}(h_n))(\hat h_n - {\hat h}_{M})
$$
converges in distribution to a centered normal law with variance $V$.

 {\textcolor{black}{Furthermore,
implicitly defining   $\delta_3(h)$ by
$G_X(h)= T_n(h)+ \delta_3(h)+ n^{-1}\|U^{1/2}(Y-r)\|^2$, if, in addition to  the 3 steps above, we show that
\begin{enumerate}
\setcounter{enumi}{3}
\item $a_n(\delta_3'({h})- \delta_2'({h}))$ converges in probability to $0$, uniformly over $H_n$
\item $\frac{G_X^{\prime\prime}(h^*)}{\BBe(T_n^{''}(h_n))}$ tends in probability to $1$ for any $h^*$ between $\hat h_n$ and $\genhGCV$,
\end{enumerate}
then
$$
a_n \BBe(T_n^{\prime\prime}(h_n))(\hat h_n - \genhGCV)
$$
converges in distribution to the same centered normal law.}}
\end{lem}
According to Lemma \ref{lemtec}, we   have to consider  {\textcolor{black}{five}} steps. We study each of them in the following   {\textcolor{black}{five}} subsections. The  {\textcolor{black}{final}} subsection concludes the proof of Theorem \ref{main}.

\subsubsection{Step 1: convergence in distribution of $a_n\delta'_2(h_n)$.}
The following proposition studies the asymptotic distribution of $a_n\delta'_2(h_n)$ for $a_n=\sqrt{\frac{n}{h_n^2}}$.
\begin{pro}\label{corocoro} Suppose that the assumptions of Theorem \ref{main} are satisfied. Then the following two assertions are equivalent.
\begin{itemize}
\item $\sqrt{\frac{n}{h_n^2}}\delta'_2(h_n)$ converges in distribution as $n$ tends to infinity to a centered  normal  law with variance  $4V$
\item $\sqrt{\frac{n}{h_n^2}} \sum_{i=1}^n\left({\tilde a}_{i,n}(h_n)u(x_i)\epsilon_i+ \sum_{j=1}^{i-1}(u(x_i)+u(x_j))b_{i,j}(h_n) \epsilon_i\epsilon_j\right)$
converges to a centered normal law with variance $V$, where
$$
{\tilde a}_{i,n}(h_n)= -C_K\frac{h_n}{n}r''(x_i).
$$
\end{itemize}
\end{pro}
\noindent {\bf{Proof of Proposition \ref{corocoro}.}}
Recall that,
for $G(u)=-uK'(u)$, for any $h\in ]0,\epsilon[$ and any $n\geq1$,
\begin{eqnarray}\label{step1}
&&
\delta_2'(h) = 2 \sum_{i=1}^na_i(h)u(x_i)\epsilon_i + 2\sum_{i=1}^n\sum_{j=1}^{i-1} (u(x_i)+u(x_j))b_{i,j}(h) \epsilon_i\epsilon_j{\nonumber}\\
&& +
2\sum_{i=1}^n b_{i,i}(h) u(x_i)(\epsilon_i^2-\BBe(\epsilon_i^2)),
\end{eqnarray}
where,
$$a_i(h) = -\frac{1}{n}\frac{\partial}{\partial h} \BBe(\hat r(x_i))= \frac{1}{n^2h^2}\sum_{j=1}^n (K-G)(\frac{x_i-x_j}{h})r(x_j), $$
$$ b_{i,j}(h) = \frac{1}{n^2h^2}K(\frac{x_i-x_j}{h})- \frac{1}{n^2h^2} G(\frac{x_i-x_j}{h}) \,\,\,{\rm and} \,\,\, b_{i,i}(h)=\frac{1}{n^2h^2}K(0).$$
We also need, for the proof of Proposition \ref{corocoro}, the following two lemmas.
\begin{lem}\label{lem6} Recall that $h_n=c n^{-1/5}$ and suppose that
$
\sum_{j=1}^{\infty}|\Cov(\epsilon_1^2,\epsilon_j^2)|<\infty,
$
then
$$
\lim_{n\rightarrow \infty} \frac{n}{h_n^2} \Var\left(\sum_{i=1}^n u(x_i) b_{i,i}(h_n) (\epsilon_i^2-\BBe(\epsilon_i^2))\right)=0.
$$
\end{lem}
{\noindent}{\bf{Proof of Lemma \ref{lem6}.}} We have,
\begin{eqnarray*}
&& \frac{n}{h_n^2} \Var\left(\sum_{i=1}^n u(x_i) b_{i,i}(h_n) (\epsilon_i^2-\BBe(\epsilon_i^2))\right)=  \frac{n}{h_n^2n^4h_n^4}K^2(0)\sum_{i=1}^n\sum_{j=1}^nu(x_i)u(x_j)\Cov(\epsilon_i^2,\epsilon_j^2) \\
&& \leq \frac{1}{n^2h_n^6}K^2(0)\|u\|_{\infty}^2 \sup_i \sum_{j=1}^{\infty}|\Cov(\epsilon_i^2,\epsilon_j^2)|.
\end{eqnarray*}
The proof of this lemma is achieved since $\lim_{n \rightarrow \infty}n^2h_n^6= \lim_{n\rightarrow \infty} n^{4/5}=\infty$.
\begin{lem}\label{term1}
Recall that $h_n=c n^{-1/5}$.  We have, noting $C_K= \int x^2K(x)dx$,
\begin{eqnarray*}
\lim_{n\rightarrow \infty}\frac{n}{h_n^2}\Var\left( \sum_{i=1}^{n}\left(a_i(h_n)+ C_K\frac{h_n}{n}r''(x_i)\right)u(x_i) \epsilon_i\right)=0.
\end{eqnarray*}
\end{lem}
\noindent {\bf Proof of Lemma \ref{term1}.}
Clearly, we have using Lemma {\bf{A.3}} of Appendix {\bf{A.2}} of the supplementary material,
$$
a_i(h_n)+ C_K\frac{h_n}{n}r''(x_i)=O(\frac{h_n^2}{n}+ \frac{1}{n^2h_n^3}),
$$
and
\begin{eqnarray*}
&& \Var\left( \sum_{i=1}^{n}(a_i(h_n)+ C_K\frac{h_n}{n}r''(x_i))u(x_i) \epsilon_i\right)\\
&&= \sum_{i=1}^n\sum_{j=1}^n(a_i(h_n)+ C_K\frac{h_n}{n}r''(x_i))(a_j(h_n)+ C_K\frac{h_n}{n}r''(x_j))u(x_i)u(x_j)\Cov(\epsilon_i,\epsilon_j)\\
&& \leq cst\,\sup_i \left((a_i(h_n)+ C_K\frac{h_n}{n}r''(x_i))u(x_i)\right)^2 n \sigma^2\\
&& = O\left(n (\frac{h_n^2}{n}+ \frac{1}{n^2h_n^3})^2\right)=O(n^{-9/5})=o(\frac{h_n^2}{n}).
\end{eqnarray*}
The proof of Lemma \ref{term1} is complete.

\noindent {\bf{{End of the proof of Proposition \ref{corocoro}.}}}  We have, using (\ref{step1}),
\begin{eqnarray*}
&& \sqrt{\frac{n}{h_n^2}}\delta'_2(h_n)= 2 \sqrt{\frac{n}{h_n^2}} \sum_{i=1}^n\left({\tilde a}_{i,n}(h_n)u(x_i)\epsilon_i+ \sum_{j=1}^{i-1}(u(x_i)+u(x_j))b_{i,j}(h_n) \epsilon_i\epsilon_j\right)\\
&& + 2 \sqrt{\frac{n}{h_n^2}} \sum_{i=1}^n\left(a_{i}(h_n)-{\tilde a}_{i,n}(h_n)\right)u(x_i)\epsilon_i+ 2 \sqrt{\frac{n}{h_n^2}} \sum_{i=1}^n b_{i,i}(h_n)u(x_i)
(\epsilon_i^2-\BBe(\epsilon_i^2)).
\end{eqnarray*}
The proof of Proposition \ref{corocoro} is complete if
  $\sqrt{\frac{{n}}{h_n^2}} \sum_{i=1}^n\left(a_{i}(h_n)-{\tilde a}_{i,n}(h_n)\right)u(x_i)\epsilon_i $  and
  $\sqrt{\frac{n}{h_n^2}} \sum_{i=1}^n b_{i,i}(h_n)u(x_i)(\epsilon_i^2-\BBe(\epsilon_i^2))$
  converge in probability to $0$ as $n$ tends to
infinity, which are satisfied due to  Lemmas  \ref{lem6} and \ref{term1}.
\subsubsection{Step 2: convergence in probability of $a_n(\delta_2'({h}_n)- \delta_2'({\hat h_n}))$}
The following proposition checks step 2 of Lemma \ref{lemtec}.
\begin{pro}\label{step2} Under the assumptions of Theorem \ref{main},
$\sqrt{\frac{n}{h_n^2}}(\delta_2'({h}_n)- \delta_2'({\hat h_n}))$ converges in probability to $0$ as $n$ tends to infinity.
\end{pro}
\noindent{\bf{Proof of Proposition \ref{step2}.}}
We have,
\begin{eqnarray*}
&& \sqrt{\frac{n}{h_n^2}}(\delta_2'({h}_n)- \delta_2'({\hat h_n}))= \sqrt{\frac{n}{h_n^2}}(h_n-{\hat h_n}) \delta_2^{\prime\prime}(h^*) = \sqrt{n}(1-\frac{{\hat h_n}}{h_n}) \delta_2^{\prime\prime}(h^*),
\end{eqnarray*}
where $h^*$ is an element of $H_n$ between $h_n$ and $\hat h_n$ and since $\hat h_n/h_n$ converges in probability to $1$ as $n$ tends to infinity (by Proposition \ref{pro1}),
we deduce that
\begin{equation}\label{prob}
\lim_{n\rightarrow \infty}\BBp(h^*\notin A_{\epsilon})=0,\,\,\, \forall\,\,\epsilon>0,
\end{equation}
where, for fixed $\epsilon>0$, $A_{\epsilon}=\{h\in H_n,\,\, |\frac{h}{h_n}-1|\leq \epsilon\}$. Now, we have for any $M>0$,
\begin{eqnarray*}
&& \BBp\left( a_n|\delta_2'({h}_n)- \delta_2'({\hat h_n})|\geq \epsilon^2\right)\\
&& \leq \BBp\left(\sqrt{n}|1-\frac{{\hat h_n}}{h_n}||\delta_2^{\prime\prime}(h^*)|\geq \epsilon^2, h^* \in A_{\epsilon}\right)+
\BBp\left(\sqrt{n}(1-\frac{{\hat h_n}}{h_n}) |\delta_2^{\prime\prime}(h^*)|\geq \epsilon^2, h^* \notin A_{\epsilon} \right)\\
&& \leq \BBp\left(\sqrt{n}|1-\frac{{\hat h_n}}{h_n}| |\delta_2^{\prime\prime}(h^*)|\BBone_{h^* \in A_{\epsilon}}\geq \epsilon^2 \right)+
\BBp\left( h^* \notin A_{\epsilon} \right)\\
&& \leq \BBp\left(\sqrt{n}|1-\frac{{\hat h_n}}{h_n}| \sup_{h\in A_{\epsilon}}|\delta_2^{\prime\prime}(h)|\geq \epsilon^2 \right)+
\BBp\left( h^* \notin A_{\epsilon} \right) 
\end{eqnarray*}
\begin{eqnarray*}
&& \leq \BBp\left(\sqrt{n}|1-\frac{{\hat h_n}}{h_n}| \sup_{h\in A_{\epsilon}}|\delta_2^{\prime\prime}(h)|\geq \epsilon^2,\,\, \sqrt{n}\sup_{h\in A_{\epsilon}}|\delta_2^{\prime\prime}(h)|\geq M \right)\\
&& +
\BBp\left(\sqrt{n}|1-\frac{{\hat h_n}}{h_n}| \sup_{h\in A_{\epsilon}}|\delta_2^{\prime\prime}(h)|\geq \epsilon^2,\,\, \sqrt{n}\sup_{h\in A_{\epsilon}}|\delta_2^{\prime\prime}(h)|< M \right)
+\BBp\left( h^* \notin A_{\epsilon} \right) \\
&& \leq  \BBp\left(\sqrt{n}\sup_{h\in A_{\epsilon}}|\delta_2^{\prime\prime}(h)|\geq M \right)+ \BBp\left(M|1-\frac{{\hat h_n}}{h_n}|\geq \epsilon^2\right)+
\BBp\left( h^* \notin A_{\epsilon} \right),
\end{eqnarray*}
which tends to $0$ by letting first $n$ tends to infinity and then $M$ tends to infinity,  due to
Proposition \ref{pro1}, (\ref{P1-2}) and (\ref{prob}).

\subsubsection{Step 3: convergence in probability of $\frac{{\rm CL}''(h^*)}{\BBe(T_n^{\prime\prime}(h_n))}$}
\begin{pro}\label{step3} Under the assumptions of Theorem \ref{main},
 $\frac{{\rm CL}''(h^*)}{\BBe(T_n^{\prime\prime}(h_n))}$ tends in probability to 1, as $n\rightarrow \infty$, for any $h^*$ between $\hat h_n$ and $\hat h_{M}$.
\end{pro}
{\noindent}{\bf{Proof of Proposition \ref{step3}.}}
We have, for any $h>0$,
$
{\rm CL}^{\prime\prime}(h)= T_n^{\prime\prime}(h)+ \delta^{\prime\prime}_2(h),
$
and
$$
\frac{{\rm CL}^{\prime\prime}(h)}{\BBe(T_n^{\prime\prime}(h))}= \frac{T_n^{\prime\prime}(h)}{\BBe(T_n^{\prime\prime}(h))}+ \frac{\delta^{\prime\prime}_2(h)}{\BBe(T_n^{\prime\prime}(h))}.
$$
Our first purpose is to prove that,
$
\sup_{h\in H_n} \frac{|\delta^{\prime\prime}_2(h)|}{|\BBe(T_n^{\prime\prime}(h))|}
$
converges to $0$, in probability, as $n$ tends to infinity. Since $
\sup_{h\in H_n} \left|\frac{D^{\prime\prime}_n(h)}{|\BBe(T_n^{\prime\prime}(h))|}-1\right|
$
converges  to $0$ as $n$ tends to infinity, with
$$
D^{\prime\prime}_n(h)=3 h^2 \int_{0}^1 u(x) r^{\prime\prime 2}(x)dx \int_{-1}^{1} t^2K(t)dt + \frac{2}{nh^3} (\int_0^1u(x)dx) \int_{-1}^1 K^2(y)dy \sigma^2,
$$
(see Lemma {\bf{A.4}} of Appendix {\bf{A.2}} of the supplementary material) it remains then to prove that, (since $\inf_{h\in H_n} nh^3|D^{\prime\prime}_n(h)|>0$),
$$
\sup_{h\in H_n} \left(nh^3|\delta^{\prime\prime}_2(h)|\right)\longrightarrow 0, \,\, {\mbox{in probability as}}\,\,\, n \rightarrow \infty,
$$
which is proved due to Lemma \ref{P1} of Subsection \ref{A.3}.
It remains to prove that
$$
\sup_{h\in H_n} \left|\frac{T_n^{\prime\prime}(h)}{\BBe(T_n^{\prime\prime}(h))}-1\right|\longrightarrow 0,\,\, \,\, {\mbox{in probability as}}\,\,\, n \rightarrow \infty,
$$
or equivalently,
\begin{eqnarray*}
\sup_{h\in H_n}{nh^3\left|T^{\prime\prime}_n(h)-\BBe\left(T^{\prime\prime}_n(h)\right)\right|}\longrightarrow 0,\,\, {\mbox{in probability as}} \,\,\, n \rightarrow \infty,
\end{eqnarray*}
which is proved due to Lemma \ref{P2} of Subsection \ref{A.3}. Consequently,
$$
\sup_{h\in H_n}\left|\frac{{\rm CL}^{\prime\prime}(h)}{\BBe(T_n^{\prime\prime}(h))}-1\right|\longrightarrow 0,\,\, {\mbox{in probability as}} \,\,\, n \rightarrow \infty,
$$
Finally,
\begin{eqnarray*}
&& \left|\frac{{\rm CL}^{\prime\prime}(h^*)}{D_n^{\prime\prime}(h_n)}-1\right|\leq cst\,\sup_{h\in H_n}\left|\frac{{\rm CL}^{\prime\prime}(h)}{D_n^{\prime\prime}(h)}-1\right|+
\left|\frac{D^{\prime\prime}_n(h^*)}{D_n^{\prime\prime}(h_n)}-1\right|.
\end{eqnarray*}
By definition of $h^*$ and by Proposition \ref{pro1}, we deduce that
$$
\left|\frac{D^{\prime\prime}_n(h^*)}{D_n^{\prime\prime}(h_n)}-1\right|\longrightarrow 0\,\, {\mbox{in probability as}}\,\,\, n \rightarrow \infty,
$$
and then
$$
\left|\frac{{\rm CL}^{\prime\prime}(h^*)}{D_n^{\prime\prime}(h_n)}-1\right|\longrightarrow 0\,\, {\mbox{in probability as}}\,\,\, n \rightarrow \infty.
$$

{\textcolor{black}{
\subsubsection{Step 4: uniform convergence in probability, over $H_n$, of $a_n  (\delta_3'({h})- \delta_2'({h}))$, where $\delta_3$ is defined in Lemma \ref{lemtec}}
The following proposition checks step 4 of Lemma \ref{lemtec} for $a_n=
n^{7/10}$ .
\begin{pro}\label{step4} Under the assumptions of Theorem \ref{main},
$n^{7/10} (\delta_3'({h})- \delta_2'({h}))$ converges in probability to $0$ {\textcolor{black}{uniformly}} over $H_n$ as $n$ tends to infinity.
\end{pro}
\noindent{\bf{Proof of Proposition \ref{step4}.}}
From the definition of $G_{\rm X}$ and the two expressions (\ref{Mallowsdep}), (\ref{delta2}) of ${\rm CL}$, we have
\begin{eqnarray*}
&& G_{\rm X}^{\prime}(h)= \frac{-K(0)}{ n h^2} \Xi_{\rm X}^\prime{\left(\frac{K(0)}{ n h}\right)} \times   \left( n^{-1} \|U^{1/2}(I-L)Y\|^2   \right)   \\
&&\,\,\,\,\,\,\,\,\,   +     \Xi_{\rm X}{\left(\frac{K(0)}{ n h}\right)}  \times  { {\rm d}  \over { {\rm d} h} }\left(n^{-1}  \|U^{1/2}(I-L)Y\|^2   \right) \\
 && =     \frac{-K(0)}{ n h^2} \Xi_{\rm X}^\prime{\left(\frac{K(0)}{ n h}\right)}   \left(  T_n(h)  +   \delta_2(h) +  n^{-1} \|U^{1/2}(Y-r)\|^2    -2  \sigma^2
 \frac{\sum_{i=1}^nu(x_i))}{ n}   {\frac{K(0)}{ n h}}  \right)     \,\, \\
 &&\,\,\,\,\,\,\,\,\, +     \Xi_{\rm X}{\left(\frac{K(0)}{ n h}\right)}    \left(  T_n^{\prime}(h)  +   \delta_2^{\prime}(h)   -   2 \sigma^2  \frac{\sum_{i=1}^nu(x_i))}{ n}    \left( {-\frac{K(0)}{ n h^2}}\right)   \right).
\end{eqnarray*}
{\textcolor{black}{Let,}}
\begin{eqnarray}
&&I_1 :=                \frac{-K(0)}{ n h^2} \Xi_{\rm X}^\prime{\left(\frac{K(0)}{ n h}\right)}   \left(  T_n(h)  +  \delta_2(h)  \right)         \\
&&I_2 := \left( \Xi_{\rm X}{\left(\frac{K(0)}{ n h}\right)}   -1 \right)   \left(  T_n^{\prime}(h)  +   \delta_2^{\prime}(h)      \right)
\end{eqnarray}
and
\begin{eqnarray}\label{I3}
I_3 :=   &&\frac{-K(0)}{ n h^2} \Xi_{\rm X}^\prime{\left(\frac{K(0)}{ n h}\right)}   \left(  n^{-1} \|U^{1/2}(Y-r)\|^2  \right)   \\
&& +
 \Xi_{\rm X}{\left(\frac{K(0)}{ n h}\right)}        \left(2  \sigma^2   \frac{\sum_{i=1}^nu(x_i))}{ n}   {\frac{K(0)}{ n h^2}}   \right).
\end{eqnarray}
From the above expression of $G_{\rm X}^{\prime}$ we can write (since, by definition $ \delta_3^{\prime}(h) = G_{\rm X}^{\prime}(h) - T_n^{\prime}(h) $)
\begin{eqnarray}\label{delata3Minusdelta2}
 \delta_3^{\prime}(h) -\delta_2^{\prime}(h) =  I_1 +I_2 +I_3  +
  2  \sigma^2   \frac{K(0)}{ n h^2} \Xi_{\rm X}^\prime{\left(\frac{K(0)}{ n h}\right)}
 \frac{\sum_{i=1}^nu(x_i))}{ n}   {\frac{K(0)}{ n h}}.
\end{eqnarray}
The least term of this decomposition is clearly $o(n^{-7/10})$ for $h$ in $H_n$ since  $\Xi_{\rm X}^\prime(t)$ is bounded near $0$. We are going to show that this {\textcolor{black}{is}} also true, in probability, uniformly over $H_n$, for the first three terms. Let us begin with $I_3$. By rewriting the factor
$ \Xi_{\rm X}{\left(\frac{K(0)}{ n h}\right)}  = (1/2) \Xi_{\rm X}^\prime{\left(\frac{K(0)}{ n h}\right)}   + \left( \Xi_{\rm X}{\left(\frac{K(0)}{ n h}\right)}    - (1/2) \Xi_{\rm X}^\prime{\left(\frac{K(0)}{ n h}\right)}  \right)$
in the second term (4.23)  of $I_3$ and invoking the required properties of  $\Xi_{\rm X}$ near $0$ (recall that  $\Xi_{\rm X}^\prime(t) = 2 + O(t)$) then
\begin{equation}\label{decompOfTermI3}
I_3 =  -  \frac{K(0)}{ n h^2} \Xi_{\rm X}^\prime{\left(\frac{K(0)}{ n h}\right)}    \ \left(  n^{-1} \|U^{1/2}(Y-r)\|^2
 -  \sigma^2   \frac{\sum_{i=1}^nu(x_i))}{ n}   \right)
      + o(n^{-2}h^{-3}).\end{equation}
%
%
%$$I_3 =  -  \frac{K(0)}{ n h^2} \Xi_{\rm X}^\prime{\left(\frac{K(0)}{ n h}\right)}    \ \left(  n^{-1} \|U^{1/2}(Y-r)\|^2
% -  \sigma^2   \frac{\sum_{i=1}^nu(x_i))}{ n}   \right)
%      + o(n^{-2}h^{-3}).$$
Now, similarly as in the proof of Lemma \ref{lem6}, it can be checked that $n^{-1} \|U^{1/2}(Y-r)\|^2 -\sigma^2 \frac{\sum_{i=1}^nu(x_i))}{ n}$ is $o_P(n^{-1/10})$,
 and thus $I_3$ is uniformly $o_P(n^{-7/10})$.
For $I_1$ and $I_2$, a such uniform rate results  from (\ref{P1-1}), (\ref{P1-1-firstDeriv}), (\ref{P2-1}), (\ref{P2-1-firstDeriv}) used  in a classical way via Markov inequality, and of the stated properties of $D_n(h)$ and $D_n^{\prime}(h)$.
}}

{\textcolor{black}{
\subsubsection{Step 5: convergence in probability of $\frac{G_{\rm X}^{\prime\prime}(h^*)}{\BBe(T_n^{\prime\prime}(h_n))}$}
We simply check  this step
%establish that this step is true
by invoking the above Step 3, and using the following result, where $\delta_3$ is defined in Lemma \ref{lemtec}, whose proof is postponed to the supplementary material  (see its Section {\bf{A.4}}) .
\begin{pro}\label{step5} Under the assumptions of Theorem \ref{main},
 $n^{2/5} ( \delta_3^{\prime\prime}(h)  - \delta_2^{\prime\prime}(h) )$ tends in probability to $0$, as $n\rightarrow \infty$, uniformly over $H_{n}$.
\end{pro}
}}

\subsubsection{End of the proof of Theorem \ref{main}}
 {\textcolor{black}{For the first part of this Theorem we}} have to check the {\textcolor{black}{first}} three items of Lemma \ref{lemtec}. We have, from Proposition {\bf{B.2}} of the supplementary material,
$$
n^{7/10} \sum_{i=1}^n\left({\tilde a}_{i,n}(h_n)u(x_i)\epsilon_i+ \sum_{j=1}^{i-1}(u(x_i)+u(x_j))b_{i,j}(h_n)
\epsilon_i\epsilon_j\right)\Longrightarrow {\cal N}(0,V).
$$
It follows from Proposition \ref{corocoro} that
$$n^{7/10}\delta'_2(h_n)\Longrightarrow {\cal N}(0,4V),$$
where $V= c^2C_K^2\sigma^2\int_0^1 u^2(x)r^{\prime\prime 2}(x)dx+ \frac{4}{c^3}\sigma^4 \int_0^1 u^2(x)dx\int_0^1 (K-G)^2(u)du$.
The two other items of Lemma \ref{lemtec} are satisfied using Propositions \ref{step2} and \ref{step3}. The proof of
{\textcolor{black}{the first part of}} Theorem \ref{main}
is complete using  Lemma \ref{lemtec} together with the fact that (see Lemma {\bf{A.4}} of the supplementary material), $\BBe(T_n^{\prime\prime}(h_n))$ is equivalent to
$$
n^{-2/5}\left(3 c^2 \int_{0}^1 u(x) r^{\prime\prime 2}(x)dx \left(\int_{-1}^{1} t^2K(t)dt\right)^2 +
\frac{2}{c^3} (\int_0^1u(x)dx) \int_{-1}^1 K^2(y)dy \sigma^2\right).
$$
{\textcolor{black}{As to the second part, arguing as above, it is now sufficient to observe  that the fourth and fifth items of Lemma \ref{lemtec} are also satisfied using Propositions \ref{step4} and \ref{step5}.}}

\section*{Acknowledgements}
{\textcolor{black}{We thank the Editor Marcus Reiss, the Associate Editor and two reviewers for their accurate and constructive comments.
This led to a significant improvement of the original manuscript, in particular in the extent of  the results, both theoretical and experimental, we present here.}}
This paper was developed in the framework of Grenoble Alpes Data Institute (ANR-15-IDEX-02).

\phantomsection
\section*{Supplement}
%\slink[doi]{COMPLETED BY THE TYPESETTER}      !!!!!!
%The supplementary material \cite{KDSsupp}   first gives, in its Section {\bf {A.1}}, the proof of Lemma \ref{ase}     (the tools used for this lemma are now classical but we give the details, there, to be self-contained). Next, it gives the proofs of the statements 
%taken for ``granted''
%(more precisely, the statements whose proofs are indicated as postponed to the supplementary material) 
%  in  the above proofs of the  Lemmas \ref{P1}, \ref{P2}, \ref{lemtec}, and \ref{term1}, the ones of Propositions \ref{step3}  and \ref{step5}, and in the final part (Section 4.3.6) of the proof of Theorem \ref{main}.

The supplementary material \cite{KDSsupp}  gives first, in its section  {\bf {A.1}}, the proof of Lemma \ref{ase}. Next, 
%in sections A.2, A.3 and A.4, 
it gives complementary proofs of Lemmas \ref{P1}, \ref{P2}, \ref{lemtec}, and \ref{term1}, of Propositions \ref{step3}  and \ref{step5}, and of Theorem \ref{main}. Moreover, some other important probabilistic properties for MDS are stated in Appendices B and C. More precisely, Appendix B gives the central limit theorem for MDS;  Appendix C gives and proves some ingredients for MDS used throughout the proofs of the main results (such as Marcinkiewicz-Zygmund type inequalities or maximal bounds for weighted sums of MDS or quadratic form of MDS).

%
%it gives complementary proofs of Lemmas 4.1, 4.2, 4.3, and 4.5, of Propositions 4.3 and 4.5, and of Theorem 3.1. Moreover, some other important probabilistic properties for MDS are stated in Appendices B and C. More precisely, Appendix B gives the central limit theorem for MDS;  Appendix C gives and proves some ingredients for MDS used throughout the proofs of the main results (such as Marcinkiewicz-Zygmund type inequalities or maximal bounds for weighted sums of MDS or quadratic form of MDS). 

%\begin{supplement}
%%\sname{Supplement A}
%%\stitle{Supplement to ``On bandwidth selection problems in nonparametric  trend estimation under martingale difference errors''}
%\stitle{Supplement}
%\slink[doi]{COMPLETED BY THE TYPESETTER}
%\sdatatype{.pdf}
%\sdescription{The supplementary material \cite{KDSsupp}   contains the proofs of Lemmas 2.1 and 4.1, 4.2, 4.3, 4.5 and the ones of Propostions 4.3 and 4.5.
%%We provide additional supporting plots that show both good and poor performance of the Hill estimator for the index of regular variation in a variety of examples.}
%\end{supplement}

\end{document}